\documentclass{amsart}

\usepackage{amssymb}
\usepackage{amsmath}
\usepackage{amsthm}
\usepackage{graphicx,psfrag}

\theoremstyle{plain}
\newtheorem{theorem}{Theorem}[section]
\newtheorem{lemma}[theorem]{Lemma}
\newtheorem{corollary}[theorem]{Corollary}
\newtheorem{proposition}[theorem]{Proposition}

\theoremstyle{definition}
\newtheorem{definition}[theorem]{Definition}
\newtheorem{notation}[theorem]{Notation}
\newtheorem*{labelling}{Labelling}
\newtheorem{claim}{Claim}

\theoremstyle{remark}
\newtheorem{remark}[theorem]{Remark}

\numberwithin{figure}{section}

\newcommand{\Int}{\mathrm{int}}

\begin{document}

\title{Heegaard surfaces and the distance of amalgamation}
\author{Tao Li}
\thanks{Partially supported by an NSF grant}

%\date{\today}
\address{Department of Mathematics \\
 Boston College \\
 Chestnut Hill, MA 02467}
\email{taoli@bc.edu}

\begin{abstract}
Let $M_1$ and $M_2$ be orientable irreducible 3--manifolds with connected boundary and suppose $\partial M_1\cong\partial M_2$.  Let $M$ be a closed 3--manifold obtained by gluing $M_1$ to $M_2$ along the boundary.  We show that if the gluing homeomorphism is sufficiently complicated, then $M$ is not homeomorphic to $S^3$ and all small-genus Heegaard splittings of $M$ are standard in a certain sense. In particular, $g(M)=g(M_1)+g(M_2)-g(\partial M_i)$, where $g(M)$ denotes the Heegaard genus of $M$.  This theorem is also true for certain manifolds with multiple boundary components.
\end{abstract}

\maketitle

\tableofcontents

\section{Introduction}\label{Sintro} 
One of the most useful ways of constructing a new 3--manifold is to glue two given 3--manifolds with boundary via a homeomorphism between their boundary surfaces.  This construction is called amalgamation.  Dehn filling and Heegaard splitting can be viewed as examples of such a construction.  In this paper, we study Heegaard splittings of 3--manifolds obtained by amalgamation.  Like Dehn filling, the 3--manifold obtained by amalgamation depends on the gluing homeomorphism.  We will show that if the gluing homeomorphism is sufficiently complicated, then the small-genus Heegaard splittings of the resulting 3--manifold are standard.

The complexity of the gluing homeomorphism is defined using the curve complex.  The curve
complex of $F$, introduced by Harvey \cite{H}, is defined as follows.  Let $F$ be a closed orientable connected surface.  The curve
complex of $F$ is the complex whose 
vertices are the isotopy classes of essential simple closed curves in
$F$.  If the genus of $F$ is at least 2, then $k+1$ vertices in the curve complex determine a $k$--simplex if they are represented by pairwise disjoint curves.  
If $F$ is a torus, then $k+1$ vertices determine a $k$--simplex if they are represented by curves that pairwise meet exactly once.  Clearly the curve complex of the torus is the same as the Farey graph.  We denote the curve complex
of $F$ by $\mathcal{C}(F)$.  For any two vertices in $\mathcal{C}(F)$,
the distance $d(x,y)$ is the minimal number of 1--simplices in a
simplicial path jointing $x$ to $y$.  To simplify notation, unless
necessary, we do not distinguish a vertex in $\mathcal{C}(F)$ from a
simple closed curve in $F$ representing this vertex.

Let $M_1$ and $M_2$ be orientable irreducible 3--manifolds with boundary. Let $F_i$ be a boundary component of $M_i$ ($i=1,2$).  In this paper, we suppose $M_i$ is not a product $F_i\times I$ and $\partial M_i-F_i$ (if not empty) is incompressible in $M_i$.  Suppose  $F_1\cong F_2\cong F$.  We can glue $M_1$ to $M_2$ via a homeomorphism $\phi:F_1\to F_2$ and obtain an orientable 3--manifold $M=M_1\cup_\phi M_2$.  We may view $M_1$ and $M_2$ as submanifolds of $M$ and $F=M_1\cap M_2$ as a closed non-peripheral surface embedded in $M$. 

\begin{definition}\label{D1}
Let $M_1$, $M_2$, $M$ and $F$ be as above.  If $F$ is compressible in $M_i$, the \emph{disk complex} of $M_i$ is the set of vertices in $\mathcal{C}(F)$ represented by curves bounding compressing disks in $M_i$.  If $M_i$ is a twisted $I$--bundle over a closed non-orientable surface, the \emph{annulus complex} of $M_i$ is the set of vertices in $\mathcal{C}(F)$ represented by boundary curves of vertical annuli in $M_i$.  If $M_i$ has incompressible boundary and $M_i$ is not a twisted $I$--bundle over a closed non-orientable surface, we fix a properly embedded essential surface $\Omega_i$ in $M_i$ with $\partial\Omega_i\cap F\ne\emptyset$ and suppose the Euler characteristic $\chi(\Omega_i)$ is maximal among all such essential surfaces.  We define $\mathcal{U}_i$ to be the set of vertices in $\mathcal{C}(F)$ as follows,
\[ \mathcal{U}_i=\left\{ \begin{array}{ll}
\text{the disk complex of $M_i$,} & \text{if $F$ is compressible in $M_i$}\\ 
\text{the annulus complex of $M_i$,} &  \text{if $M_i$ is a twisted $I$--bundle}\\
\text{vertices represented by components of $\partial\Omega_i\cap F$,} & \text{otherwise.}
\end{array} \right.\]
We define the distance of the amalgamation to be $d(M)=d(\mathcal{U}_1,\mathcal{U}_2)$ in the curve complex $\mathcal{C}(F)$.
\end{definition}

Note that the surface $\Omega_i$ in Definition~\ref{D1} is not unique, but we will show in section~\ref{Ssmall} that, if $M_i$ has incompressible boundary and is not an $I$--bundle, then the diameter of the set of vertices in $\mathcal{C}(F)$ represented by boundary curves of such essential surfaces is bounded.  Thus any different choice of $\Omega_i$ only changes $d(M)$ by an explicit small number.  
If both $M_1$ and $M_2$ are handlebodies or more generally if $F$ is compressible in both $M_1$ and $M_2$, then $d(M)$ is the same as the Hempel distance, see \cite{H, ST1}.   Schleimer informed the author that, similar to the disk complex, the annulus complex of a twisted $I$--bundle is also quasi-convex in $\mathcal{C}(F)$.  So $d(M)$ is arbitrarily large if the gluing map $\phi$ is a sufficiently high power of a pseudo-Anosov map.  Like the Hempel distance, $d(M)$ also provides a natural complexity measure for a one-sided Heegaard splitting, i.e., a decomposition of $M$ into a handlebody and a twisted $I$--bundle.

If one prefers, the following is a roughly equivalent way of defining $d(M)$ which does not involve a choice of $\Omega_i$.  Let $k_i$ be the maximal Euler characteristic of essential orientable surfaces properly embedded in $M_i$ and with at least one boundary component in $F$ (we consider compressing disks as essential surfaces).  Let $\mathcal{U}_i$ be the set of vertices in $\mathcal{C}(F)$ represented by boundary curves of such essential surfaces whose Euler characteristic is $k_i$. Then one can define $d(M)=d(\mathcal{U}_1,\mathcal{U}_2)$.  If $k_i=1$, then $\mathcal{U}_i$ is the disk complex of $M_i$. If $M_i$ has incompressible boundary and is not a twisted $I$--bundle, then by section~\ref{Ssmall}, the diameter of $\mathcal{U}_i$ is bounded by a number depending only on $k_i$.

\begin{theorem}\label{TS3}
Let $M=M_1\cup_F M_2$ be as above.  Then there is a number $K$ depending on $M_1$ and $M_2$ such that if $d(M)\ge K$ then 
\begin{enumerate}
\item $M$ is irreducible and $\partial$-irreducible, and 
\item $M$ is not homeomorphic to $S^3$.
\end{enumerate}
\end{theorem}

Similar to Dehn surgery, one can perform a surgery on a graph in $S^3$.  An immediate corollary of Theorem~\ref{TS3} is that given a graph $\Gamma$ in $S^3$, if one performs a complicated surgery on $\Gamma$, i.e., gluing back a handlebody to $S^3-N(\Gamma)$ via a high-distance map, then the resulting closed 3-manifold is irreducible and cannot be $S^3$.

\begin{definition}\label{Dmiddle}
Let $N$ be a compression body and $F$ a closed separating surface properly embedded in $N$. $F$ cuts $N$ into two submanifolds $N_1$ and $N_2$ and we may view $F$ as a boundary component of each $N_i$.  Suppose $F$ is not a 2--sphere and $\partial_+N\subset\partial N_1$.  We say $F$ is a \emph{middle surface} in $N$ if both $N_1$ and $N_2$ are compression bodies, $\partial_+N=\partial_+N_1$, $\partial_-N_2\subseteq\partial_-N$ and $F=\partial_+N_2\subseteq\partial_-N_1$.  Note that if one views a compression body as a manifold obtained by adding 2--handles and 3--handles to $\partial_+N\times I$ on the same side, then a middle surface is a middle level of this process.  In particular, one can find a handle structure of $N$ such that $N_1$ is a compression body obtained by adding a subset of the 2-- and 3--handles to $\partial_+N\times I$, and after adding the remaining 2-- and 3--handles along $F$, we obtain the whole of $N$.  Note that unless $F$ is parallel to a component of $\partial_- N$, $F$ is incompressible in $N_1$ but compressible in $N_2$.
\end{definition}

Next we consider the untelescoping of a Heegaard splitting, see \cite{ST, S2}. Let $M=V\cup_S W$ be an irreducible Heegaard splitting.  We may view the compression body $V$ as the manifold obtained by attaching 1--handles to either a product neighborhood of $\partial_-V$ or to a 0--handle; and view $W$ as the manifold obtained by attaching 2--handles and possibly a 3--handle to a product neighborhood of $S=\partial_+W$.  So a Heegaard splitting gives a natural handle-decomposition of $M$.  The untelescoping of the Heegaard splitting is a rearrangement of the order in which these handles are attached.  This rearrangement gives a decomposition of $M$ into submanifolds $N_1,\dots, N_m$ along incompressible surfaces, and each $N_i$ inherits a strongly irreducible Heegaard splitting from a subset of the original 1-- and 2--handles, see \cite{ST} for details. The decomposition is often called a generalized Heegaard splitting.  We summarize this as the following theorem.    Note that by part (1) of Theorem~\ref{TS3}, if the gluing map is sufficiently complicated, then $M$ has incompressible boundary.  So in this paper we only consider the case that $\partial M$ (if not empty) is incompressible, though untelescoping is also defined for manifolds with compressible boundary.  

\begin{theorem}[Scharlemann-Thompson \cite{ST}]\label{TST}
Let $M$ be an irreducible and orientable 3--manifold with incompressible boundary.  Let $S$ be an unstabilized Heegaard surface of $M$.  Then the untelescoping of the Heegaard splitting described above gives a decomposition of $M$ as follows, see Figure~\ref{Funtel} for a picture.
\begin{enumerate}
\item $M=N_0\cup_{F_1}N_1\cup_{F_2}\dots\cup_{F_m}N_m$, where each $F_i$ is incompressible in $M$.  
\item Each $N_i=A_i\cup_{P_i}B_i$, where each $A_i$ and $B_i$ is a union of compression bodies with $\partial_+A_i=P_i=\partial_+B_i$ and $\partial_-A_i=F_i=\partial_-B_{i-1}$.  
\item Each component of $P_i$ is a strongly irreducible Heegaard surface of a component of $N_i$.
\item no component of $A_i$ and $B_i$ is a trivial compression body (i.e. a product).
\item the genus $g(F_i)<g(S)$ and $g(P_i)\le g(S)$ for each $i$.
\end{enumerate}
\end{theorem}

Let $F$ be a closed connected surface embedded in $M$.  We say $F$ is a \emph{canonical surface with respect to the untelescoping} if $F$ is parallel to a middle surface in a component of $A_i$ or $B_i$, for some $i$.  Note that a component of $P_i$ is a (trivial) middle surface for both $A_i$ and $B_i$, and any component of $F_i$ is a middle surface for both $B_{i-1}$ and $A_i$.  The main theorem of the paper is:

\begin{figure} 
\begin{center}
\psfrag{dots}{$\dots\dots$}
\psfrag{N0}{$N_0$}
\psfrag{N1}{$N_1$}
\psfrag{Nm}{$N_m$}
\psfrag{F1}{$F_1$}
\psfrag{F2}{$F_2$}
\psfrag{Fm}{$F_m$}
\psfrag{P0}{$P_0$}
\psfrag{P1}{$P_1$}
\psfrag{Pm}{$P_m$}
\psfrag{A0}{$A_0$}
\psfrag{A1}{$A_1$}
\psfrag{A'}{$A_m$}
\psfrag{B0}{$B_0$}
\psfrag{B1}{$B_1$}
\psfrag{Bm}{$B_m$}
\includegraphics[width=3in]{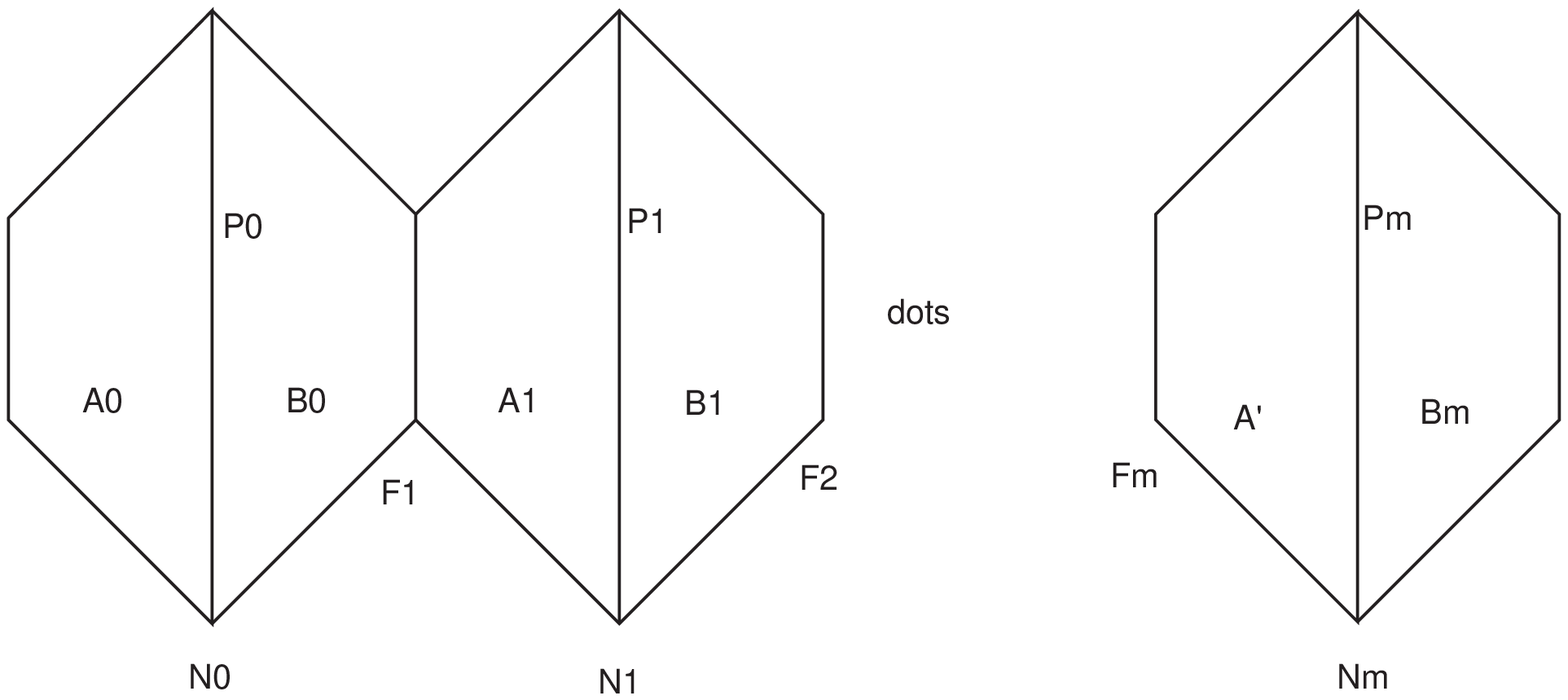}
\caption{}\label{Funtel}
\end{center}
\end{figure}

\begin{theorem}\label{Tmain}
Let $M=M_1\cup_F M_2$ be as above. In particular, suppose $M_i$ is not a product $F\times I$ and suppose $\partial M_i-F$ (if not empty) is incompressible in $M_i$.   Then for any integer $g$, there is a number $K$ depending on $M_1$, $M_2$ and $g$, such that if $d(M)>K$, then for any unstabilized Heegaard surface $S$ of $M$ with $g(S)\le g$, $F$ is isotopic to a canonical surface with respect to any untelescoping of $S$.
\end{theorem}

A corollary of Theorem~\ref{Tmain} is a formula for the Heegaard genus.

\begin{corollary}\label{Cgenus}
Let $M=M_1\cup_F M_2$ be as in Theorem~\ref{Tmain}.  Then there is a number $K$ depending on $M_1$ and $M_2$ such that if $d(M)>K$, $g(M)=g(M_1)+g(M_2)-g(F)$.
\end{corollary}

It follows from the proof that the number $K$ in Theorem~\ref{Tmain} and Corollary~\ref{Cgenus} can be chosen to be an explicit quadratic function of $g$ and $\chi(\Omega_i)$, where $\Omega_i$ is as in Definition~\ref{D1}.  The bound $K$ depends on several distance estimates in various places in the paper.  Lemma~\ref{Ldist} is the only place where the bound is quadratic and all other estimates are linear functions.

If both $M_1$ and $M_2$ are simple, i.e., irreducible, $\partial$-irreducible, atoroidal and anannular, then Theorem~\ref{Tmain} is a generalization of a theorem of Lackenby \cite{La} and is proved in \cite{So, L4}.  Note that the complexity measure in \cite{L4} is defined using boundary curves of normal surfaces, see \cite{L1} for a relation between Heegaard surfaces and normal surfaces.   

Our motivation for the main theorem is to study the Heegaard genus of closed 3--manifolds.  The following are special cases of Theorem~\ref{Tmain}.  

\begin{corollary}\label{Cmain}
Let $M_1$ and $M_2$ be orientable irreducible 3--manifolds with connected boundary and suppose $\partial M_1\cong\partial M_2$.  Let $M$ be a closed 3--manifold obtained by gluing $M_1$ to $M_2$ along the boundary.  Then for any integer $g$, there is a number $K$ depending on $M_1$, $M_2$ and $g$, such that if $d(M)>K$, then for any unstabilized Heegaard surface $S$ of $M$ with $g(S)\le g$, $F$ is isotopic to a canonical surface with respect to any untelescoping of $S$.  Moreover, $g(M)=g(M_1)+g(M_2)-g(F)$.
\end{corollary}

\begin{corollary}\label{Chandle}
Let $M$, $M_1$ and $M_2$ be as above.  If $M_1$ is a handlebody, then there is a number $K$ depending on $M_2$ and $g$ such that if $d(M)>K$, any Heegaard surface of $M$ with genus at most $g$ is isotopic to a Heegaard surface of $M_2$.  In particular, $g(M)=g(M_2)$.
\end{corollary}

Corollary~\ref{Chandle} says that if the gluing map is complicated, then there is no new small-genus Heegaard surface in the resulting 3--manifold $M$.  In particular, if $M_1$ is a solid torus and $M_2$ is a knot manifold, Corollary~\ref{Chandle} gives a weaker version of a known result on Heegaard structure and Dehn filling, see \cite{MR, MS, RSed}.

Theorem~\ref{Tmain} and its proof give useful ways of constructing new 3--manifolds with certain control on the Heegaard genus.  This may shed light on constructing counterexamples (or examples) for the rank conjecture, which asserts that for a closed (hyperbolic) 3--manifold, the rank of its fundamental group equals its Heegaard genus.  For example, if one can construct an example of 3--manifold $N$ with connected boundary whose rank is smaller than its Heegaard genus, then by Corollary~\ref{Chandle}, one can obtain a closed 3--manifold $\hat{N}$ by capping off $\partial N$ using a handlebody and via a sufficiently complicated gluing map, such that $rank(\hat{N})<g(\hat{N})$.  
Very recently, using hyperbolic JSJ pieces, the author has constructed examples of closed 3--manifolds with rank smaller than genus \cite{L7}.  These are the first such examples having hyperbolic JSJ pieces and a main tool in the construction is Theorem~\ref{Tmain}.  It is conceivable that this method may be generalized to give a hyperbolic counterexample to the rank conjecture.

The proof of Theorem~\ref{Tmain} is also used by Bachman in \cite{B} to study stabilization of Heegaard splittings. 

In the proof of the main theorem, we study how the amalgamation surface $F$ intersects the incompressible and strongly irreducible surfaces in the untelescoping of a Heegaard splitting.  In particular, we will show that if the amalgamation distance $d(M)$ is sufficiently large, then any small-genus closed incompressible surfaces in $M$ can be isotoped disjoint from $F$.  So by our definition of $d(M)$, it is natural to study the distance between $\partial\Omega_i$ and the boundary curves of incompressible or strongly irreducible surfaces in $M_i$ ($i=1,2$).

A basic observation in the proof is that the diameter in the curve complex of the vertices represented by the boundary curves of certain small-genus surfaces is bounded.  More precisely, let $N$ be a 3--manifold with incompressible boundary (which is not an $I$-bundle) and let $F$ be a boundary component of $N$.  We consider the set of surfaces in $N$ that are either essential or strongly irreducible and $\partial$-strongly irreducible (see Definition~\ref{D33}) with bounded genus and bounded number of components in $\partial N-F$.  The observation is that the diameter in the curve complex $\mathcal{C}(F)$ of the boundary curves of such surfaces is bounded.  This is proved in section~\ref{Ssmall}.

In section~\ref{SI}, we use this observation to prove Theorem~\ref{Tmain} in the case that $F$ is incompressible in both $M_1$ and $M_2$.  The case that $F$ is compressible in both $M_1$ and $M_2$ basically follows from a theorem of Scharlemann and Tomova \cite{ST1} and this is discussed in section~\ref{SII}.  The most difficult case in Theorem~\ref{Tmain} is that $F$ is compressible on one side but incompressible on the other side. The last two sections are devoted to this case.

In section~\ref{Scase3}, we show that if the distance of the amalgamation $d(M)$ is large, then $F$ can be isotoped disjoint from all the incompressible surfaces (i.e. the $F_i$'s in Theorem~\ref{TST}) in the untelescoping.  In section~\ref{Scase3}, we also discuss the case that $F$ is disjoint from all the strongly irreducible Heegaard surfaces $P_i$'s, and prove that in a certain generic situation, $F$ is isotopic to a middle surface of a compression body in the Heegaard splitting of $N_i$ (see Theorem~\ref{TST}).  In section~\ref{Slast}, we study how $F$ intersects the sweepout of the strongly irreducible Heegaard splitting of $N_i$ in the untelescoping and finish the proof of the last case.

I would like to thank Saul Schleimer for a helpful conversation on the annulus complex of a twisted $I$--bundle and thank the referee for many suggestions and corrections.

\section{A genus calculation}

\begin{notation}
Throughout this paper, we denote the interior of $X$ by $\Int(X)$, the closure of $X$ (under the path-metric) by $\overline{X}$, and the number of components of $X$ by $|X|$ for any space $X$.
\end{notation}

We first show that Corollary~\ref{Cgenus} follows from Theorem~\ref{Tmain}.  For simplicity, we suppose $M_i$ and $M_2$ have only one boundary component, i.e., $M$ is a closed 3--manifold.  If $M$ has boundary, by Theorem~\ref{TS3},  $\partial M$ is incompressible in $M$ and one may cap off each component of $\partial M$ by a handlebody and calculate the Heegaard genus same as the case that $M$ is closed.

Suppose $M$ is closed. 
By Theorem~\ref{TS3}, we may assume $M=M_1\cup_FM_2$ is irreducible and is not $S^3$.  Let $S$ be an unstabilized Heegaard surface of $M$. Let $M=N_0\cup_{F_1}N_1\cup_{F_2}\dots\cup_{F_m}N_m$ and $N_i=A_i\cup_{P_i}B_i$ be the decompositions in an untelescoping of the Heegaard splitting, see Theorem~\ref{TST} and Figure~\ref{Funtel}.  
As in \cite{ST, S2}, one can rearrange the handle structure determined by the Heegaard splitting along $S$ so that the sub-collection of 1-- and 2--handles which occur in $N_i$ determine the Heegaard splitting $N_i=A_i\cup_{P_i}B_i$.  

Suppose $S$ is a minimal genus Heegaard surface of $M$ and let $g$ be its genus.  Suppose $F$ is canonical with respect to the untelescoping of $S$ as above.  Without loss of generality, we may suppose $F$ lies in the compression body $B_j$ between $P_j$ and $F_{j+1}$ in the untelescoping, see Figure~\ref{Funtel}.  We may assume $B_j$ is connected.  By the definition of middle surface, $F$ separates $B_j$ into two compression bodies and we can choose a handle structure for $B_j$ so that the 2--handles in the two compression bodies are exactly the 2--handles for $B_j$.  Next we count the handles in $M_1$ and $M_2$.

Let $a_i$, $b_i$, $c_i$ and $d_i$ ($i=1,2$) be the numbers of 0--, 1--, 2--, and 3--handles in $M_i$ respectively in the handle decomposition determined by the Heegaard surface $S$ as above.  The total number of 0--handles is $a_1+a_2$ and the total number of 1--handles is $b_1+b_2$.  So the Heegaard genus  $g=(b_1+b_2)-(a_1+a_2)+1=(c_1+c_2)-(d_1+d_2)+1$.  

Since $F$ and $M_1$ are connected, as in \cite{ST}, one can rearrange the 0-- and 1--handles in $M_1$ to form a connected handlebody and obtain a Heegaard splitting of $M_1$ with genus $g_1=b_1-a_1+1$.  Hence $g(M_1)\le b_1-a_1+1$.  Similarly, one can rearrange the 2-- and 3--handles in $M_2$ to form a handlebody and obtain a Heegaard splitting of $M_2$ with genus $g_2=c_2-d_2+1$. Hence $g(M_2)\le c_2-d_2+1$.  Moreover, an easy calculation of the Euler characteristic of $M_1$ yields $g(F)=1-a_1+b_1-c_1+d_1$.  Therefore, 
$g(M_1)+g(M_2)-g(F)\le (b_1-a_1+1)+(c_2-d_2+1)-(1-a_1+b_1-c_1+d_1)=(c_1+c_2)-(d_1+d_2)+1=g=g(M).$

Given two minimal-genus Heegaard splittings of $M_1$ and $M_2$, the amalgamation of the two splittings yields a Heegaard splitting of $M$ with genus $g(M_1)+g(M_2)-g(F)$, see \cite{La, L4, Sch} for more detailed description. This means that $g(M)\le g(M_1)+g(M_2)-g(F)$.  So the equality $g(M)= g(M_1)+g(M_2)-g(F)$ holds.

\section{Intersection of small surfaces}\label{Ssmall}

In this section, we prove several lemmas on the intersection of certain small-genus surfaces.  These lemmas will be used in the later sections.

Throughout this section, we fix an orientable irreducible compact connected 3-manifold $N$ with incompressible boundary.  We also fix a component of $\partial N$ and denote it by $F$.

\begin{definition}\label{Dann}
We define the annulus complex $\mathcal{A}_N(F)$ to be the subcomplex of $\mathcal{C}(F)$ consisting of vertices represented by boundary curves of essential annuli in $N$. Note that we only consider those essential annuli with at least one boundary component in $F$.
\end{definition}

The following lemma is also proved in \cite{L5}.

\begin{lemma}\label{Lannulus}
Suppose $N$ is not an $I$--bundle. Then the diameter of the annulus complex $\mathcal{A}_N(F)$ in $\mathcal{C}(F)$ is at most 2.
\end{lemma}
\begin{proof}
Let $J$ be an $I$--bundle in $N$ with its horizontal boundary $\partial_hJ$ in $\partial N$ and its vertical boundary consisting of essential annuli properly embedded in $N$.  Suppose $\partial_hJ\cap F\ne\emptyset$.  Note that if $N$ contains an essential annulus $A$ with at least one boundary component in $F$, then a small neighborhood of $A$ is such an $I$--bundle.  We may suppose $J$ is maximal up to isotopy.  This is basically from the theory of characteristic submanifolds, see \cite{Ja}.  

As $N$ is not an $I$--bundle, $J\ne N$.  By our assumption, $J\cap F\subset\partial_hJ$ and any component of $\partial (J\cap F)$ is a boundary component of an essential annulus in $N$.
Let $A'$ be a vertical boundary component of $J$.  So $A'$ is an essential annulus in $N$ and $\partial A'\cap F\subset \partial (J\cap F)$.  
Let $A$ be any other essential annulus in $N$ with at least one boundary component in $F$ and we consider $A\cap A'$.  Since $A$ and $A'$ are both essential annuli, no component of $A\cap A'$ can be essential in one annulus but trivial in the other annulus.  If $A\cap A'$ contains a closed curve that is trivial in both annuli, then there is such a curve $c$ that is innermost in $A'$ and bounding disks $\Delta\subset A$ and $\Delta'\subset A'$.  Since $c$ is innermost in $A'$, $\Int(\Delta')\cap A=\emptyset$ and $\Delta\cup\Delta'$ is an embedded $S^2$ in $N$.  Since $N$ is irreducible,  $\Delta\cup\Delta'$ must bound a 3-ball.  Hence we can perform an isotopy on $A$, pushing $\Delta$ across the 3-ball and eliminate the intersection curve $c$.  So after isotopy, we may assume $A\cap A'$ contains no trivial closed curve. 
 If $A\cap A'$ contains an arc that is trivial in both annuli, then there is such an arc $\alpha$ that is outermost in $A'$.  Since $\alpha$ is trivial in both $A$ and $A'$, $\alpha$ and subarcs of $\partial A$ and $\partial A'$ bound bigon disks $d$ and $d'$ in $A$ and $A'$ respectively.  Moreover, since $\alpha$ is outermost in $A'$ and $A\cap A'$ contains no trivial closed curve, $\Int(d')\cap A=\emptyset$ and $d\cup d'$ is a disk properly embedded in $N$.  Since $\partial N$ is incompressible, the disk $d\cup d'$ must be $\partial$-parallel in $N$.  Hence an isotopy on $A$ that pushes $d$ across the 3-ball bounded by $d\cup d'$ and $\partial N$ can eliminate the intersection arc $\alpha$.   Thus after some isotopies as above, every arc or closed curve in $A\cap A'$ is essential in both annuli and this means that either $\partial A\cap\partial A'=\emptyset$ or $A\cap A'$ consists of arcs vertical in both $A$ and $A'$.  

If $\partial A\cap\partial A'\ne\emptyset$ after isotopy, then the union of a small neighborhood of $J\cup A$ and possibly some 3--balls yields a larger $I$--bundle contradicting the assumption that $J$ is maximal, see \cite[Section 2]{L0} for a more detailed argument.  So $\partial A\cap\partial A'=\emptyset$ after isotopy.  This means that, for any component $\gamma$ of $\partial (J\cap F)$, $d(\gamma, \partial A\cap F)\le 1$ and the lemma holds.
\end{proof}

\begin{definition}\label{D33}
Let $N$ be an orientable irreducible compact connected 3-manifold with incompressible boundary as above. Let $Q$ be a surface properly embedded in $N$ and suppose $Q$ is not a disk or 2-sphere.
We say $Q$ is \emph{essential} if it is incompressible and $\partial$-incompressible.  A properly embedded disk in $N$ is essential if its boundary is an essential curve in $\partial N$, and a 2-sphere in $N$ is essential if it does not bound a 3-ball in $N$.  Since $N$ is irreducible and $\partial N$ is incompressible, $N$ contains no essential disk or 2-sphere. 
Let $P$ be a properly embedded separating surface in $N$ and we allow $P$ to be disconnected.  Suppose the surface $P$ decomposes $N$ into two submanifolds $X$ and $Y$, where $X$ and $Y$ are on different sides of $P$ (note that $X$ and $Y$ may be disconnected). 
We say $P$ is \emph{strongly irreducible} if $P$ has compressing disks on both sides, and each compressing disk in $X$ meets each compressing disk in $Y$.  We say $P$ is \emph{$\partial$-strongly irreducible} if 
\begin{enumerate}
\item every compressing and $\partial$-compressing disk in $X$ meets every compressing and $\partial$-compressing disk in $Y$, and
\item there is at least one compressing or $\partial$-compressing disk on each side of $P$.
\end{enumerate}
\end{definition}
If $P$ is strongly irreducible, then $\partial P$ consists of curves essential in $\partial N$. To see this, suppose a component of $\partial P$ is trivial in $\partial N$. Then an innermost such component bounds a disk in $\partial N$ that is disjoint from every compressing disk on the other side of $P$.  This contradicts that $P$ is strongly irreducible.

Let $P$ be a strongly irreducible and $\partial$-strongly irreducible surface in $N$ and  $\partial P\ne\emptyset$. Let $X$ and $Y$ be the closure of the two submanifolds of $N-P$ on different sides of $P$ as in Definition~\ref{D33}.  Since $P$ is compressible on both sides, we may compress $P$ in both $X$ and $Y$.  Let $P^X$ and $P^Y$ be the possibly disconnected surfaces obtained by maximally compressing $P$ in $X$ and $Y$ respectively and removing all possible 2--sphere components.  Some components of $P^X$ and $P^Y$ may be closed surfaces.  Let $P^X_\partial $ (resp. $P^Y_\partial $) be the union of the components of $P^X$ (resp. $P^Y$) with boundary.  

For any $\partial$-parallel surface $R$ in $N$, we denote by $\pi(R)$ the subsurface of $\partial N$ that is bounded by $\partial R$ and isotopic to $R$ relative to $\partial R$.  We say a collection of pairwise disjoint $\partial$-parallel surfaces $R_1,\dots, R_m$ in $N$ are \emph{non-nested} if $\pi(R_1),\dots,\pi(R_m)$ are pairwise disjoint in $\partial N$.

\begin{lemma}\label{Lside}
Let $P$, $P^X$, $P^Y$, $P^X_\partial $ and $P^Y_\partial $ be as above.  Then 
\begin{enumerate}
\item[(a)] $P^X$ and $P^Y$ are incompressible in $N$,
\item[(b)] a component of $P^X_\partial $ is either $\partial$-parallel or can be changed into an essential surface after some $\partial$-compressions in $X$ and deleting any resulting $\partial$-parallel components, 
\item[(c)] the $\partial$-parallel components of $P^X_\partial $ are non-nested in $X$.
\end{enumerate}  
\end{lemma}
\begin{proof}
Since $P$ is strongly irreducible, part (a) of the lemma follows from \cite[Lemma 5.5]{S1}.  Our task is to prove parts (b) and (c).

As $P$ is separating, we call the two sides of $P$ plus and minus sides and suppose $P^X$ is on the plus side and $P^Y$ is on the minus side.  Moreover, any surface obtained by compression or $\partial$-compression on $P$ inherits plus and minus sides.
So $P^X$ is a surface obtained by compressing $P$ on the plus side. 
 
Let $Q$ be either $P^X$ or a surface obtained from $P^X$ by some $\partial$-compressions on the plus side (i.e. in $X$), and let $Q'$ be a component of $Q$ with boundary.  By part (a), $P^X$ is incompressible, hence $Q$ and $Q'$ are incompressible in $N$.

Next prove that $Q'$ is $\partial$-incompressible on the minus side.  The main reason for this is that $P$ is $\partial$-strongly irreducible and the proof is similar to \cite[Lemma 5.5]{S1}. 

Suppose $Q'$ is $\partial$-compressible on the minus side and let $D$ be a $\partial$-compressing disk for $Q'$ on the minus side.  By viewing $P$ as a surface obtained from $Q$ by adding some tubes and possibly some half tubes on the minus side, we may view the arc $\alpha=\partial D\cap Q'$ as an arc in $P$ and view $D$ as a disk transverse to $P$ with $\partial D=\alpha\cup\beta$, $\alpha\subset P$, $\beta\subset\partial N$, and $\partial\alpha=\partial\beta$. Moreover, a neighborhood of $\alpha$ in $D$ lies on the minus side of $P$.  Note that if $\Int(D)\cap P=\emptyset$, then $D$ is a $\partial$-compressing disk for $P$ on the minus side disjoint from a compressing disk of $P$ on the plus side, contradicting that $P$ is $\partial$-strongly irreducible. So $\Int(D)\cap P\ne\emptyset$. 
Let $\gamma_1,\dots,\gamma_n$ be the closed curves in $\Int(D)\cap P$ and let $\alpha_1,\dots,\alpha_k$ be the arcs in $(D-\alpha)\cap P$ with $\partial\alpha_i\subset\beta$ for each $i$.   After isotopy, we may assume the $\gamma_i$'s and $\alpha_i$'s are essential curves and arcs in $P$. 

Since $P$ is strongly irreducible, it follows from the proof of Scharlemann's no-nesting lemma \cite[Lemma 2.2]{S} (also see \cite[Lemma 5.5]{S1}), after some isotopy (one can also use the isotopy described below), we may assume the closed curves $\gamma_i$'s are not nested in $D$.  Let $\delta_1,\dots,\delta_n$ be the subdisks of $D$ bounded by $\gamma_1,\dots,\gamma_n$ respectively.  Each $\alpha_i$ and a subarc of $\beta$ bound a subdisk $D_i$ of $D$.  Since $P$ can be viewed as the surface obtained from $Q$ by adding tubes and half tubes on the minus side corresponding to the compressions and $\partial$-compressions on the plus side, we may assume that (1) each $\delta_i$ is a compressing disk for $P$ in $X$ (i.e. on the plus side), and (2) if $\Int(D_i)\cap P=\emptyset$, $D_i$ is a $\partial$-compressing disk for $P$ in $X$ (i.e. on the plus side).

If some $\delta_i$'s lie inside some $D_j$, since the $\delta_i$'s are non-nested in $D$, there must be a disk $D_j$ such that $\Int(D_j)\cap P\ne\emptyset$ but those disks $D_i$'s and $\delta_i$'s that lie inside $D_j$ are pairwise disjoint (i.e.~non-nested in $D_j$).  This assumption implies that a small neighborhood of $\alpha_j$ in $D_j$ lies on the minus side of $P$ (since those $\delta_i$'s and $D_i$'s in $D_j$ are in $X$).  Thus, after replacing $D$ by this disk $D_j$ in our argument if necessary, we may assume all the disks $D_i$'s and $\delta_i$'s are non-nested in $D$.  Furthermore, we may assume $|\Int(D)\cap P|$, the number of components of $\Int(D)\cap P$, is minimal among all such disks $D$.  Note that the isotopies above eliminating  nested closed curves and trivial intersection curves all reduce $|\Int(D)\cap P|$.   So each $\delta_i$ is a compressing disk for $P$ on the plus side and each $D_i$ is a $\partial$-compressing disk for $P$ on the plus side.

Since $P$ is compressible on both sides, $P$ has a compressing disk $D'$ on the minus side.  Since $P$ is strongly irreducible and $\partial$-strongly irreducible, $\partial D'\cap\partial\delta_i\ne\emptyset$ and $\partial D'\cap\partial D_i\ne\emptyset$ for each $i$.  Thus $D'\cap D\ne\emptyset$.  We may assume $D'$ is transverse to $D$.

After some isotopies, we may also assume $D'\cap D$ does not contain any closed curve.  We may assume $|D\cap D'|$ is minimal up to isotopy.  Let $\kappa$ be an arc in $D'\cap D$ that is outermost in $D'$, i.e., $\kappa$ and a subarc of $\partial D'$ bound a subdisk $\Delta$ of $D'$ and $\Int(\Delta)\cap D=\emptyset$.  We have the following 8 cases to consider.

\emph{Case (1)}. If $\kappa$ is an arc connecting two different circles $\gamma_i$ and $\gamma_j$ in $D$, then a simple isotopy that pushes $D$ across $\Delta$ will merge $\gamma_i$ and $\gamma_j$ into one closed curve. This contradicts the assumption that $|\Int(D)\cap P|$ is minimal. 

\emph{Case (2)}. If $\kappa$ is an arc connecting a circle $\gamma_i$ and an arc $\alpha_j$, then the same isotopy above  merges $\gamma_i$ and $\alpha_j$ into a single arc. This again contradicts that  $|\Int(D)\cap P|$ is minimal.

\emph{Case (3)}. The third case is that $\partial\kappa$ lies in the same circle $\gamma_i=\partial\delta_i$, as shown in Figure~\ref{Fcirc}(a).  After the same isotopy pushing $D$ across $\Delta$, the disk $\delta_i$ becomes an annulus $A\subset D$, see Figure~\ref{Fcirc}(b).  Moreover, $A$ is properly embedded in $X$ on the plus side of $P$ (since $\delta_i\subset X$).  

We denote the two circles of $\partial A$ by $c_1$ and $c_2$, as shown in Figure~\ref{Fcirc}(b).  Let $d_i$ be the disk bounded by $c_i$ in $D$ and suppose $d_1\subset d_2$ and $d_2-\Int(d_1)=A$.  If $c_i$ is a trivial curve in $P$, then a simple isotopy on $P$ and $D$ can eliminate $c_i$, and $|D\cap D'|$ is reduced after all these operations while $|\Int(D)\cap P|$ is either reduced or unchanged.  So we may assume both $c_1$ and $c_2$ are essential curves in $P$.  If $\Int(d_1)\cap P=\emptyset$, then $d_1$ is a compressing disk in $Y$ (since $A\subset X$) and $d_1$ can be isotoped disjoint from $\delta_i$, a contradiction to the hypothesis that $P$ is strongly irreducible.   Thus we may assume $\Int(d_1)\cap P\ne\emptyset$.

Since the circles $\gamma_i$'s are non-nested, $c_1$ and those $\gamma_j$'s in $\Int(d_1)$ bound a planar surface $R\subset d_1$ and $R$ is properly embedded in $Y$, see Figure~\ref{Fcirc}(b). By a theorem of Scharlemann \cite[Theorem 2.1 and Lemma 2.2]{S} (also see \cite[Lemma 5.5]{S1}), one can perform an isotopy to eliminate the nested circles in $d_2$.  In fact, it follows from \cite[Theorem 2.1 and Lemma 2.2]{S} and \cite[Lemma 5.5]{S1} that $R$ must be $\partial$-parallel in $Y$.  This contradicts the minimality assumption on $|\Int(D)\cap P|$.  

\begin{figure} 
\begin{center}
\psfrag{(a)}{(a)}
\psfrag{(b)}{(b)}
\psfrag{(c)}{(c)}
\psfrag{(d)}{(d)}
\psfrag{a}{$a$}
\psfrag{b}{$b$}
\psfrag{gi}{$\gamma_i$}
\psfrag{ai}{$\alpha_i$}
\psfrag{aj}{$\alpha_j$}
\psfrag{k}{$\kappa$}
\psfrag{A}{$A$}
\psfrag{R}{$R$}
\psfrag{c1}{$c_1$}
\psfrag{c2}{$c_2$}
\psfrag{D_a}{$D_a$}
\psfrag{Db}{$D_b$}
\includegraphics[width=4in]{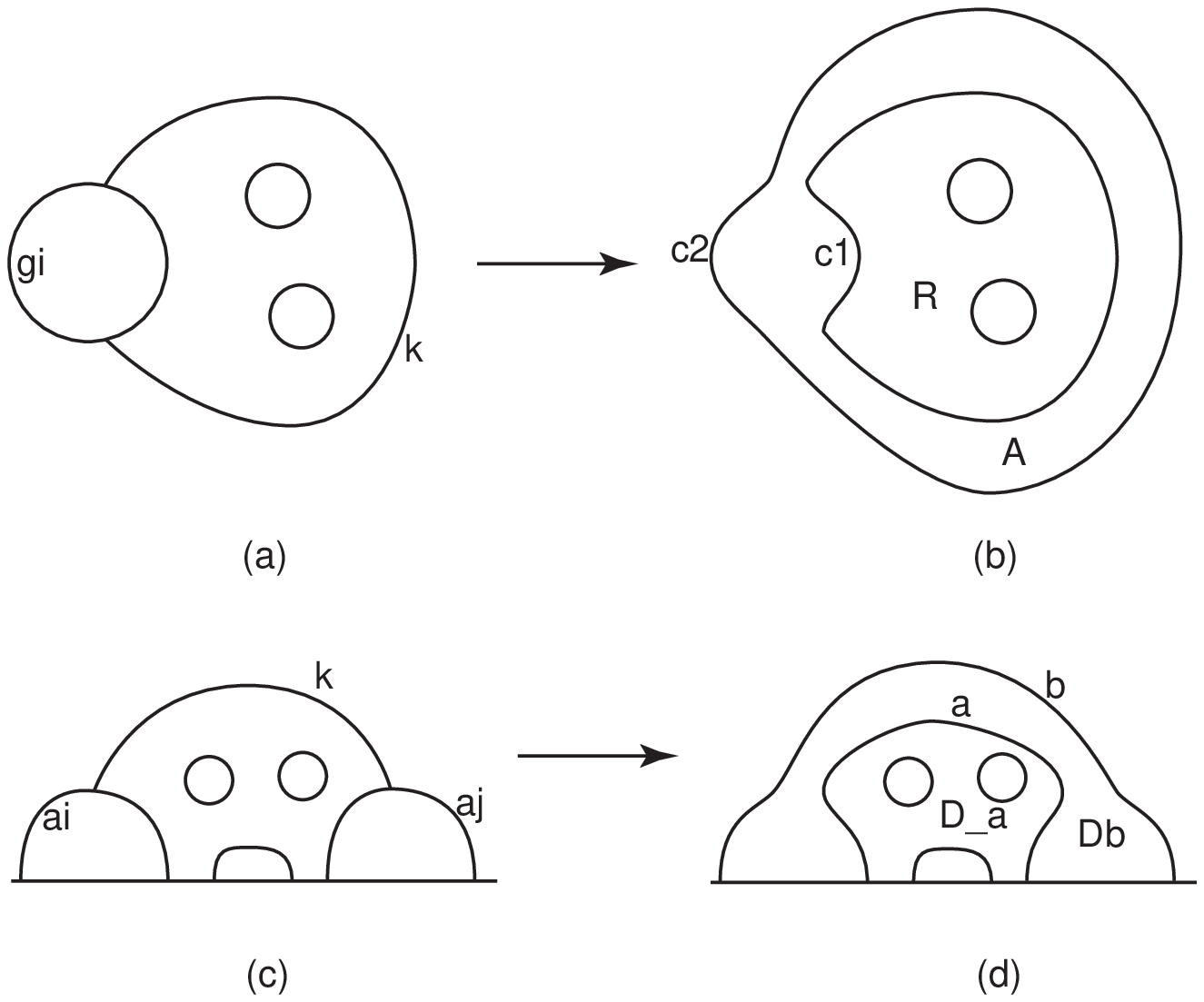}
\caption{}\label{Fcirc}
\end{center}
\end{figure}

\emph{Case (4)}. Now we consider the case that $\partial\kappa$ lies in the same arc $\alpha_i$.  After the isotopy pushing $D$ across $\Delta$ as above, $\alpha_i$ splits into an arc $\alpha_i'$ and a circle $c$.  Let $d$ be the subdisk of $D$ bounded by $c$.  By applying the same arguments for $c_1$ and $d_1$ in case (3) to $c$ and $d$, we get a contradiction to either the minimality of $|\Int(D)\cap P|$ or the assumption that $P$ is strongly irreducible.

\emph{Case (5)}. In this case we suppose $\kappa$ is an arc connecting two different arcs $\alpha_i$ and $\alpha_j$, see Figure~\ref{Fcirc}(c).  After the isotopy pushing $D$ across $\Delta$, $\alpha_i$ and $\alpha_j$ become a pair of arcs $a$ and $b$ with $\partial a\cup\partial b=\partial\alpha_i\cup\partial\alpha_j$, as shown in Figure~\ref{Fcirc}(d).  Let $D_a$ and $D_b$ be subdisks of $D$ cut off by $a$ and $b$ respectively.  By our construction, $D_a$ and $D_b$ are nested. Suppose $D_a\subset D_b$.  If $a$ is a trivial arc in $P$, then a simple isotopy on $D$ can remove the intersection arc $a$ and lead to a contradiction to the minimality of $|\Int(D)\cap P|$.  Suppose $a$ is essential in $P$.  If $\Int(D_a)\cap P=\emptyset$, then $D_a$ is a $\partial$-compressing disk for $P$ on the minus side, since $D_i$ and $D_j$ are on the plus side.  Moreover, we can perturb $D_a$ to be disjoint from $D_i$ and $D_j$, and this contradicts that $P$ is $\partial$-strongly irreducible.  Thus $\Int(D_a)\cap P\ne\emptyset$.  Since $D_i$ and $D_j$ lie on the plus side of $P$, a neighborhood of $a$ in $D_a$ is on the minus side of $P$.  Since $|\Int(D_a)\cap P|<|\Int(D)\cap P|$, this again contradicts the minimality assumption on $|\Int(D)\cap P|$. 

\emph{Case (6)}.  The final 3 cases deal with the situation that an endpoint of $\kappa$ lies in $\alpha=\partial D\cap P$.  We first suppose $\kappa$ connects $\alpha$ and a circle $\gamma_i$ in $D$. After pushing $D$ across $\Delta$ as above, $\alpha$ and $\gamma_i$ merge into a boundary arc $\alpha'$ of a new disk $E$, where $\partial\alpha'=\partial\alpha$, $\beta=\partial D\cap\partial N=\partial E\cap\partial N$ and $\partial E=\alpha'\cup\beta$.  Note that after compressing $P$ into $P^X$, the circle $\gamma_i$ can be viewed as a trivial circle in $Q'$ bounding a disk (in $Q'$) corresponding to the compressing disk, hence $\alpha'$ can be viewed as an arc in $Q'$ isotopic (in $Q'$) to $\alpha$.  However, $|\Int(E)\cap P|<|\Int(D)\cap P|$ and we have a contradiction.

\emph{Case (7)}. If $\partial\kappa\subset\alpha$, then similar to Case (4), after pushing $D$ across $\Delta$, $D$ splits into two disks $E'$ and $E''$, and $\alpha$ splits into a closed curve $\gamma'=\partial E'$ and an arc $\alpha'\subset\partial E''$. Similar to Case (6), we may view $\gamma'$ and $\alpha'$ as curves in $Q'$.  The closed curve $\gamma'$ must be trivial in $Q'$, since $Q'$ is incompressible.  Hence $\alpha'$ is isotopic in $Q'$ to $\alpha$ and the new disk $E''$ can be viewed as a $\partial$-compressing disk for $Q'$ isotopic to $D$.  However, $\kappa$ is eliminated and $|D\cap D'|$ is reduced after the isotopy while $|\Int(E'')\cap P|\le|\Int(D)\cap P|$.  This contradicts our minimality assumption on $|D\cap D'|$.

\emph{Case (8)}. If $\kappa$ connects $\alpha$ to an arc $\alpha_i$ in $D$, then after pushing $D$ across $\Delta$, similar to Case (5), $\alpha$ and $\alpha_i$ merge into a pair of new arcs $\alpha_a$ and $\alpha_b$ which are boundary arcs of two new disks $E_a$ and $E_b$ respectively ($\partial\alpha\cup\partial\alpha_i=\partial\alpha_a\cup\partial\alpha_b$, $E_a$ and $E_b$ correspond to the two components of $D-D_i-\kappa$).  Recall that the arc $\alpha_i$ bounds a $\partial$-compressing disk $D_i$ for $P$ and $P^X$.  By our construction of $Q'$, after the $\partial$-compressions on $P^X$ that we performed to get $Q$, we may view $\alpha_i$ as a $\partial$-parallel arc in $Q'$ that cuts off a disk in $Q'$ corresponding to the $\partial$-compressing disk $D_i$.  Since $\alpha$ is an essential arc in $Q'$ and $\alpha_i$ can be viewed as a trivial arc in $Q'$, at least one of $\alpha_a$ and $\alpha_b$ is an essential arc in $Q'$ bounding a $\partial$-compressing ($E_a$ or $E_b$) for $Q'$.  After replacing $D$ by a new disk $E_a$ or $E_b$ above, we get a contradiction to the minimality assumption of $|\Int(D)\cap P|$. 

Therefore $Q'$ must be $\partial$-incompressible on the minus side.  
In particular each component of $P^X_\partial $ is $\partial$-incompressible on the minus side.
Since $P^X_\partial $ is incompressible by part (a), a component of $P^X_\partial $ is either $\partial$-parallel in $X$ or can be changed into an essential surface after some $\partial$-compressions on the plus side and deleting any resulting $\partial$-parallel components.  So part (b) of the lemma holds.   If the $\partial$-parallel components of $P^X_\partial $ are nested in $X$ (i.e.~the two product regions in $X$ bounded by two $\partial$-parallel components are nested), then this means that after some $\partial$-compressions on $P^X_\partial $ on the plus side, the resulting surface becomes $\partial$-compressible on the minus side, a contradiction to the conclusion above.  Thus part (c) holds.
\end{proof}

\begin{lemma}\label{Lessential}
Let $P$ and $Q$ be properly embedded orientable surfaces in $N$ with at least one boundary component in $F$.  Suppose $Q$ is essential and $P$ is strongly irreducible and $\partial$-strongly irreducible.  Note that by the definition of strongly irreducible surface, $P$ is separating. Then either
\begin{enumerate}
\item $(\partial P\cap F)\cap(\partial Q\cap F)=\emptyset$ after isotopy, or 
\item after some compressions and $\partial$-compressions on the same side of $P$, one can obtain an essential surface with a boundary component in $F$, or 
\item after some isotopy, $P\cap Q$ is essential in both $P$ and $Q$ and $|\partial P\cap\partial Q|$ is minimal among curves isotopic to $\partial P$ and $\partial Q$.
\end{enumerate}
\end{lemma}
\begin{proof} Suppose part (1) of the lemma is not true. So
we may assume that $(\partial P\cap F)\cap(\partial Q\cap F)\ne\emptyset$ and $|\partial P\cap\partial Q|$ is minimal among curves isotopic to $\partial P$ and $\partial Q$ in $\partial N$.

Let $X$ and $Y$ be the closure of the two submanifolds of $N-P$ as in Definition~\ref{D33}.  Let $P^X$ and $P^Y$ be the possibly disconnected surfaces obtained by maximally compressing $P$ in $X$ and $Y$ respectively and removing all possible 2--sphere components.   Let $P^X_\partial $ and $P^Y_\partial $ be the unions of the components of $P^X$ and $P^Y$ with boundary respectively.  Since $P$ has at a boundary component in $F$, $P^X_\partial \cap F\ne\emptyset$ and $P^Y_\partial \cap F\ne\emptyset$.

By Lemma~\ref{Lside}, a component of $P^X_\partial $ is either $\partial$-parallel in $X$ or can be changed to an essential surface by some $\partial$-compressions on the $X$-side.  Furthermore, by Lemma~\ref{Lside}, the $\partial$-parallel components of $P^X_\partial $ are non-nested.  Thus either part (2) of Lemma~\ref{Lessential} holds or the components of $P^X_\partial$ and $P^Y_\partial $ incident to $F$ are $\partial$-parallel and non-nested in $X$ and $Y$ respectively.   
Let $P^X_F$ and $P^Y_F$ be the components of $P^X_\partial$ and $P^Y_\partial $ respectively whose boundary lie in $F$.  Suppose part (2) of Lemma~\ref{Lessential} is not true, then as above, $P^X_F$ and $P^Y_F$ are $\partial$-parallel and non-nested in $X$ and $Y$ respectively.  Note that by our hypotheses and assumptions, $P^X_F\ne\emptyset$ and $P^Y_F\ne\emptyset$.

Let $N_P$ be the submanifold of $N$ between $P^X$ and $P^Y$ and we may assume $P$ is properly embedded in $N_P$.   By the construction of $P^X$ and $P^Y$, there are graphs $G_X\subset X\cap N_P$ and $G_Y\subset Y\cap N_P$, which correspond to the compressions on $P$ in $X$ and $Y$ respectively, such that $N_P-(P^X\cup G_X\cup P^Y\cup G_Y)$ is a product $P\times(0,1)$.  Let $\Sigma_X=P^X\cup G_X$ and $\Sigma_Y=P^Y\cup G_Y$.  We may view this as a sweepout $H: P\times (I,\partial I)\to (N_P, \Sigma_X\cup\Sigma_Y)$, where $I=[0,1]$ and $H|_{P\times(0,1)}$ is an embedding.  We denote $H(P\times\{a\})$ by $P_a$ for any $a\in I$.  We may assume $P_0=\Sigma_X$ and $P_1=\Sigma_Y$ and each $P_a$ ($a\ne 0,1$) is isotopic to $P$. 

Since $P^X_F$ and $P^Y_F$ are $\partial$-parallel and non-nested, and by our hypothesis that $Q$ has a boundary component in $F$, we may assume $Q\cap P^X_F$ and $Q\cap P^Y_F$ consist of non-nested $\partial$-parallel arcs in $Q$.  Since we have assumed at the beginning that $|\partial P\cap\partial Q|$ is minimal among curves isotopic to $\partial P$ and $\partial Q$, we may assume the arcs in $Q\cap P^X_F$ and $Q\cap P^Y_F$ are essential in $P^X_F$ and $P^Y_F$ respectively. 
Moreover, after isotopy, we may assume that for each $t\in(0,1)$, $Q$ is transverse to $P_t$ except for at most one center or saddle tangency. We call $t\in (0,1)$ a regular level if $Q$ is transverse to $P_t$, otherwise, we call $t$ a singular level. We may assume there are only finitely many singular levels.

Let $X_t$ and $Y_t$ ($t\in(0,1)$) be the closure of the two submanifolds of $N-P_t$ corresponding to $X$ and $Y$ respectively.   
For each regular level $t$, we label it $X$ (resp.~$Y$) if either a closed curve in $Q\cap P_t$ bounds a compressing disk for $P_t$ in $X_t$ (resp.~$Y_t$), or an arc in $Q\cap P_t$ bounds a $\partial$-compressing disk for $P_t$ in $X_t$ (resp.~$Y_t$).  

Recall that we have assumed that $Q\cap P^X_F$ and $Q\cap P^Y_F$ consist of non-nested $\partial$-parallel  arcs in $Q$ (since $P^X_F$ and $P^Y_F$ are $\partial$-parallel in $X$ and $Y$ respectively) and the arcs in $Q\cap P^X_F$ and $Q\cap P^Y_F$ are essential in $P^X_F$ and $P^Y_F$ respectively.  This means that, for any sufficiently small $\epsilon>0$, $\epsilon$ is labelled $X$ and $1-\epsilon$ is labelled $Y$.

Since $P$ is strongly irreducible and $\partial$-strongly irreducible, no regular level $t$ is labelled both $X$ and $Y$.

Since $Q$ is an essential surface, no arc or curve in $P_t\cap Q$ can be trivial in $P_t$ but nontrivial in $Q$.  For any regular level $t$, if a closed curve or an arc is trivial in both $P_t$ and $Q$, since $N$ is irreducible and $\partial$-irreducible, this curve or arc can be eliminated by an isotopy.  So for any regular level $t$, if a closed curve (resp. an arc) in $P_t\cap Q$ is trivial in $Q$ but essential in $P_t$, then we can find an innermost (resp. outermost) such curve (resp. arc) in $Q$ that bounds a compressing disk (resp.~$\partial$-compressing disk) for $P_t$ and hence $t$ is labelled $X$ or $Y$.  Thus if a level $t$ has no label, after some isotopies removing curves and arcs trivial in both $P_t$ and $Q$, $P_t\cap Q$ consists of  curves and arcs essential in both $P_t$ and $Q$ and part (3) of the lemma holds.  So to prove the lemma, it remains to consider the case that every regular level $t$ is labelled. 

Since $\epsilon$ is labelled $X$ and $1-\epsilon$ is labelled $Y$ for small $\epsilon>0$, the conclusions above imply that there must be a singular level $s\in (0,1)$ such that $s-\epsilon$ is labelled $X$ but $s+\epsilon$ is labelled $Y$ for sufficiently small $\epsilon>0$.  Moreover, $P_s\cap Q$ contains a single saddle tangency.

Let $\Theta$ be the graph component of $P_s\cap Q$ containing the saddle tangency and let $N(\Theta)$ be the closure of a small regular neighborhood of $\Theta$ in $Q$.  Since $P_s$, $Q$ and $N$ are all orientable and $P_s$ is separating in $N$, every component of $P_{s\pm\epsilon}\cap  Q$ is isotopic in $Q$ to either a component of $\partial N(\Theta)$ or a component of $P_s\cap Q-\Theta$.  
Since $s-\epsilon$ is labelled $X$ and $s+\epsilon$ is labelled $Y$, there are arcs or closed curves $\gamma_X$ and $\gamma_Y$ in $P_{s\pm\epsilon}\cap  Q$ bounding compressing or $\partial$-compressing disks in $X_{s-\epsilon}$ and $Y_{s+\epsilon}$ respectively.  As above, since $P_s$ is separating and in particular $P_s$ separates $P_{s-\epsilon}$ and $P_{s+\epsilon}$ in $N$, $\gamma_X$ and $\gamma_Y$ correspond to disjoint curves in $\partial N(\Theta)\cup(P_s\cap Q-\Theta)$.  Let $Q\times J$ be a small product neighborhood of $Q$ in $N$, where $J$ is a closed interval. Let $Q^+$ and $Q^-$ be the two components of $Q\times\partial J$.  Note that by the configuration near a saddle tangency, the intersection patterns of $P_s\cap Q^\pm$ and $P_{s\pm\epsilon}\cap Q$ are the same.  In particular, curves in $\partial N(\Theta)\cup(P_s\cap Q-\Theta)\subset Q$ are isotopic in $Q\times J$ to disjoint curves in $P_s\cap (Q\times\partial J)$.  This means that there are two disjoint arcs or closed curves $\gamma_X'$ and $\gamma_Y'$  in $P_s\cap (Q\times\partial J)$ corresponding to $\gamma_X$ and $\gamma_Y$ above, such that $\gamma_X'$ and $\gamma_Y'$ bound compressing or $\partial$-compressing disks for $P_s$ in $X_s$ and $Y_s$ respectively. 
This contradicts that $P_s$ is strongly irreducible and $\partial$-strongly irreducible. 
\end{proof}

For any integers $g$ and $b$, let $C_{g,b}$ be the collection of orientable surfaces properly embedded in $N$, such that any $P\in C_{g,b}$ has at least one boundary component in $F$, $\partial P$ is essential in $\partial N$, $g(P)\le g$ and $|\partial P-F|\le b$.  Note that surfaces in $C_{g,b}$ need not to be essential and there is no restriction on the number of components of $\partial P\cap F$.  

\begin{lemma}\label{Lcompress}
Let $P$ be a surface in $C_{g,b}$.  Let $P=P_0, P_1,\dots, P_k$ be surfaces in $N$ such that each $P_i$ is obtained by performing a $\partial$-compression on $P_{i-1}$.  Suppose $\partial P_i$ is essential in $\partial N$ for each $i$.  Then the distance $d(\partial P\cap F,\partial P_k\cap F)\le\max\{1,\ 4g+2b-2\}$ in $\mathcal{C}(F)$. 
\end{lemma}
\begin{proof} Since $\partial N$ is incompressible in $N$ and $\partial P_i$ is essential in $\partial N$, $\chi(P_i)\le 0$ for each $i$. 
Let $b_F$ be the number of components of $\partial P\cap F$.  Since $P\in C_{g,b}$, the total number of boundary components of $P$ is at most $b_F+b$.  By our hypotheses, the total number of $\partial$-compressions is at most $-\chi(P)$, so $k\le-\chi(P)\le 2g-2+b+b_F$.

Let $D_i$ be the $\partial$-compressing disk for $P_i$ such that $P_{i+1}$ is obtained by the $\partial$-compression along $D_i$.  Let $\partial D_i=\alpha_i\cup\beta_i$ with $\alpha_i\subset P_i$ and $\beta_i\subset\partial N$.  Since we are only concerned about how the curves change in $\mathcal{C}(F)$, we may assume $\beta_i\subset F$ for all $i$.  Clearly, for any components $\gamma_i$ and $\gamma_{i+1}$ of $\partial P_i\cap F$ and $\partial P_{i+1}\cap F$ respectively, $d(\gamma_i,\gamma_{i+1})\le 1$.  Note that if $\partial\alpha_i\cap\gamma_i=\emptyset$, then the $\partial$-compression does not change $\gamma_i$ and $\gamma_i$ can be viewed as a component of $\partial P_{i+1}$. 

We may view each $\alpha_i$ above as an arc properly embedded in $P=P_0$.  As above, we may assume the endpoints of these $\alpha_i$'s all lie in $\partial P\cap F$.  We have $k$ such arcs $\alpha_i$ and $|\partial P\cap F|=b_F$.  Hence there is a component $\gamma$ of $\partial P\cap F$ that contains at most $2k/b_F$ endpoints of these arcs $\alpha_i$'s.  Since $k\le 2g-2+b+b_F$, if $2g-2+b\ge 0$, then the number of those endpoints in $\gamma$ is at most 
$$\frac{2k}{b_F}\le\frac{4g-4+2b+2b_F}{b_F}=\frac{4g-4+2b}{b_F}+2\le (4g-4+2b)+2=4g+2b-2.$$ 
If $2g-2+b<0$, we have 
$\frac{2k}{b_F}\le\frac{4g-4+2b+2b_F}{b_F}=\frac{4g-4+2b}{b_F}+2<2$, which means that the number of those endpoints in $\gamma$ is at most 1.
 
This means that at most $\max\{1,\ 4g+2b-2\}$ $\partial$-compressions occur at the curve $\gamma$.  Since each $\partial$-compression changes a curve by at most distance one in $\mathcal{C}(F)$, $d(\gamma,\partial P_k\cap F)\le\max\{1,\ 4g+2b-2\}$ and hence the lemma holds.  Note that this bound is not sharp and one can easily reduce the bound by a more delicate argument.
\end{proof}

\begin{lemma}\label{Ldist} 
Suppose $N$ is not an $I$--bundle. 
Let $P$ and $Q$ be surfaces in  $C_{g,b}$.  Suppose $Q$ is essential and suppose $P$ is either essential or strongly irreducible and $\partial$-strongly irreducible.  Then there exists a number $K'$ that depends only on $g$ and $b$, such that the distance $d(\partial P\cap F,\partial Q\cap F)\le K'$ in $\mathcal{C}(F)$.  Moreover, $K'$ can be chosen to be an explicit quadratic function of $g$ and $b$.
\end{lemma}
\begin{proof}
Before we proceed, we would like to mention a well-known result on the relation between the intersection number of two curves and their distance in the curve complex. Let $\alpha$ and $\beta$ be two essential simple closed curves in $F$ and let $k$ be the minimal number of intersection points of $\alpha\cap\beta$ up to isotopy on $\alpha$ and $\beta$ in $F$.  If the genus of $F$ is at least two, by \cite[Lemma 2.1]{MM1} (see \cite[Lemma 2.1]{H} for a better bound), the distance $d(\alpha,\beta)\le 2k+1\le 2|\alpha\cap\beta|+1$.  The proof of \cite[Lemma 2.1]{MM1} also works in the case that $F$ is a torus.  For completeness, we include the proof.  Suppose $F$ is a torus and $\alpha$ and $\beta$ realize the minimal intersection number $k$.  If $k=1$, then by the definition of the curve complex for torus, $d(\alpha,\beta)=1$.  If $k\ge 2$, we fix two points $x$ and $y$ of $\alpha\cap\beta$ adjacent in $\alpha$ and let $\alpha'$ be the subarc of $\alpha$ between $x$ and $y$ with $\Int(\alpha')\cap\beta=\emptyset$.  We can do surgery at $x$ and $y$, replacing a segment of $\beta$ between $x$ and $y$ by $\alpha'$.  This operation produces a simple closed curve $\beta_1$ in $F$ whose minimal intersection number with $\alpha$ is at most $k-1$.  Moreover, since $F$ is a torus and $\alpha\cap\beta$ realizes the minimal intersection number, the intersection points of $\alpha\cap\beta$ all have the same sign.  This implies that (1) $\beta_1$ must be nontrivial in $F$ and (2) the minimal intersection number of $\beta$ and $\beta_1$ is one, which means $d(\beta_1,\beta)=1$. Thus $d(\alpha,\beta)\le d(\alpha,\beta_1)+d(\beta_1,\beta)=d(\alpha,\beta_1)+1$.  So by inductively applying the argument above, we have $d(\alpha,\beta)\le k\le |\alpha\cap\beta|$ if $F$ is a torus.  Therefore, for any essential simple closed curves $\alpha$ and $\beta$ in $F$, we have $d(\alpha,\beta)\le 2|\alpha\cap\beta|+1$ if $g(F)\ge 2$ and $d(\alpha,\beta)\le |\alpha\cap\beta|$ if $F$ is a torus. 

We have two cases to consider.

\vspace{8pt}
\noindent
\emph{\textbf{Case A}}. $N$ contains an essential annulus with at least one boundary component in $F$.
\vspace{8pt}

Let $X$ be the set of essential annuli and M\"{o}bius bands in $N$ such that for each $A$ in $X$, (1) at least one boundary component of $A$ lies in $F$ and (2) after isotopy, $P\cap A$ is essential in both $P$ and $A$.  Note that if $P$ is an essential surface, then every essential annulus or M\"{o}bius band satisfies property (2) above. We first consider the subcase that $X\ne\emptyset$.

\vspace{8pt}
\noindent
\emph{Claim 1}.  For any annulus or a M\"{o}bius band $A_X$ in $X$, $d(\partial P\cap F, \partial A_X\cap F)\le K_1$ for some constant $K_1$ which can be chosen to be a linear function of $g$ and $b$.
\begin{proof}[Proof of the Claim 1]
Let $A$ be the annulus or M\"{o}bius band in $X$ with the properties that $P\cap A$ consists of arcs essential in both $P$ and $A$, and that $|A\cap P|$ is minimal among all such essential annuli and M\"{o}bius bands.

Note that for any essential M\"{o}bius band, the boundary of its small neighborhood in $N$ gives an essential annulus disjoint from the M\"{o}bius band.  
By Lemma~\ref{Lannulus} the diameter of the annulus complex is at most 2.  So if $d(\partial P\cap F, \partial A\cap F)\le K_1$, then for any other annulus or M\"{o}bius band $A_X$ in $X$, $d(\partial P\cap F, \partial A_X\cap F)\le K_1+3$.  Thus to prove the claim, it suffices to show that $d(\partial P\cap F, \partial A\cap F)\le K_1$ for this particular annulus or M\"{o}bius band $A$. 

If $P$ is an annulus, then by Lemma~\ref{Lannulus}, $d(\partial P\cap F, \partial A\cap F)\le 2$.  So we may assume $\chi(P)<0$. 

Suppose $d(\partial P\cap F,\partial A\cap F)\ge 2$, then every component of $\partial P\cap F$ intersects every component of $\partial A\cap F$. 
Let $\omega=\min\{|\alpha\cap\beta|:$ $\alpha$ is a component of $\partial P\cap F$ and $\beta$ is a component of $\partial A\cap F$\}.  By the result we mentioned at the beginning of the proof, $d(\partial P\cap F, \partial A\cap F)\le 2\omega+1$ if $g(F)\ge 2$ \cite[Lemma 2.1]{MM1} and $d(\partial P\cap F, \partial A\cap F)\le \omega$ if $F$ is a torus.  So to prove the claim, it suffices to show that $\omega$ is bounded from above by a linear function of $g$ and $b$.

Let $g_P$ be the genus of $P$ and let $b_F$ and $b_P$ be the numbers of components of $\partial P\cap F$ and $\partial P-F$ respectively.  By the definition of $C_{g,b}$, $g_P\le g$ and $b_P\le b$.  Since each component of $\partial P\cap F$ intersects every component of $\partial A\cap F$ in at least $\omega$ points, the number of arcs of $P\cap A$ with an endpoint in $F$ is at least $\omega b_F/2$.  
Since $\chi(P)<0$, there are at most $6g_P+3(b_F+b_P)-6\le 3b_F+6g+3b-6$ pairwise nonparallel arcs in $P$.  Thus if $\omega>\frac{2(3b_F+6g+3b-6)}{b_F}$, at least two arcs of $P\cap A$ are parallel in $P$ and each of the two arcs has at least one endpoint in $F$.  
Similar to the proof of Lemma~\ref{Lcompress}, if $6g+3b-6\ge 0$, $\frac{2(3b_F+6g+3b-6)}{b_F}= 6+\frac{12g+6b-12}{b_F}\le 12g+6b-6$, and if $6g+3b-6<0$, $\frac{2(3b_F+6g+3b-6)}{b_F}<6$.  Thus if $\omega>\max\{5,\ 12g+6b-6\}$ then there are 2 arcs in $P\cap A$, denoted by $\alpha$ and $\beta$, such that $\alpha$ and $\beta$ are parallel in $P$, $\partial\alpha\cap F\ne\emptyset$ and $\partial\beta\cap F\ne\emptyset$.

Let $R_P$ and $R_A$ be rectangles bounded by $\alpha$ and $\beta$ in $P$ and $A$ respectively. We may choose $\alpha$ and $\beta$ so that $\Int(R_P)\cap A=\emptyset$.  Thus $A'=R_P\cup R_A$ is an embedded annulus or M\"{o}bius band.  Next we show that $A'$ is an essential annulus or M\"{o}bius band.

Since $N$ is irreducible, if $A'$ is a M\"{o}bius band, then $\partial A'$ must be essential in $\partial N$ (otherwise the union of $A'$ and a disk bounded by $\partial A'$ is an embedded projective plane in $N$).  If $A'$ is an annulus and a component of $\partial A'$ is trivial in $\partial N$, then since $\partial N$ is incompressible, both components of $\partial A'$ must be trivial in $\partial N$ and this means that $\Sigma_A=(A-R_A)\cup R_P$ is an annulus isotopic to $A$ ($\Sigma_A$ is embedded because $\Int(R_P)\cap A=\emptyset$).  Moreover, after a slight perturbation on $\Sigma_A$, $\Sigma_A$ becomes transverse to $P$ with $\Sigma_A\cap P$ essential in both $\Sigma_A$ and $P$ and $|\Sigma_A\cap P|<|A\cap P|$.  This contradicts our choice of $A$.  So if $A'$ is an annulus, $\partial A'$ consists of essential curves in $\partial N$.  Since $\partial N$ is incompressible, the argument above says that no matter whether $A'$ is an annulus or a M\"{o}bius band, $A'$ is incompressible.  
Since $A$ is essential and $\alpha$ is an essential arc in $A$, $\alpha$ must be an essential arc in $N$, which implies that $A'$ is $\partial$-incompressible.  Therefore $A'$ must be an essential annulus or M\"{o}bius band in $N$.  Moreover, after a small perturbation, $A'\cap P$ has fewer components than $A\cap P$, a contradiction to the minimality assumption on $|A\cap P|$.  This means that $\omega\le\max\{5, 12g+6b-6\}$ and the claim holds.
\end{proof}

\noindent
\emph{Claim 2}.  For any essential annulus $A_N$ in $N$ with at least one boundary component in $F$, $d(\partial P\cap F, \partial A_N\cap F)\le K_2$ for some constant $K_2$ which can be chosen to be a linear function of $g$ and $b$.
\begin{proof}[Proof of the Claim 2]
By Lemma~\ref{Lannulus} the diameter of the annulus complex is at most 2.  So Claim 2 immediately follows from Claim 1 if $X\ne\emptyset$. 

If $P$ is essential in $N$, for any essential annulus $A_N$ in the claim, after isotopy, $P\cap A_N$ is essential in both $P$ and $A_N$. So $A_N\in X$ and $X\ne\emptyset$ in this subcase and Claim 2 follows from Claim 1.

If $P$ is strongly irreducible and $\partial$-strongly irreducible, by Lemma~\ref{Lessential}, either $(\partial A_N\cap F)\cap(\partial P\cap F)=\emptyset$, or one can obtain an essential surface $P'$ (with $\partial P'\cap F\ne\emptyset$) by compressing and $\partial$-compressing $P$ on the same side, or after isotopy $A_N\cap P$ consists of essential arcs in both $A_N$ and $P$.  If $(\partial A_N\cap F)\cap(\partial P\cap F)=\emptyset$, then $d(\partial A_N\cap F,\partial P\cap F)\le 1$ and the claim holds.  If $A_N\cap P$ consists of essential arcs in both $A_N$ and $P$, then $A_N\in X$ and Claim 2 follows from Claim 1.  Thus it remains to consider the subcase that one can obtain an essential surface $P'$ (with $\partial P'\cap F\ne\emptyset$) by compressing and $\partial$-compressing $P$ as in Lemma~\ref{Lessential}.  Since a compression does not change the boundary curve,  the essential surface $P'$ is obtained by $\partial$-compressing a surface whose boundary is the same as $\partial P$.  By Lemma~\ref{Lcompress}, $d(\partial P'\cap F,\partial P\cap F)\le\max\{1,\ 4g+2b-2\}$.  Since $P'$ is an essential surface, we may assume $P'\cap A_N$ is essential in both $P'$ and $A_N$.  By Claim 1,  $d(\partial P'\cap F, \partial A_N\cap F)\le K_1$.  As $d(\partial P'\cap F,\partial P\cap F)\le\max\{1,\ 4g+2b-2\}$, we have 
$$d(\partial P\cap F, \partial A_N\cap F)\le d(\partial P\cap F,\partial P'\cap F)+1+d(\partial P'\cap F, \partial A_N\cap F)\le K_2,$$
where $K_2=\max\{1,\ 4g+2b-2\}+K_1+1$.
\end{proof}

The argument above also implies that $d(\partial Q\cap F,\partial A_N\cap F)\le K_2$, since $Q$ is an essential surface.  Therefore, if $N$ contains an essential annulus with at least one boundary component in $F$, $d(\partial P\cap F,\partial Q\cap F)\le 2K_2+1$ and the lemma holds in Case A.

\vspace{8pt}
\noindent
\emph{\textbf{Case B}}. $N$ contains no essential annulus with a boundary component in $F$.
\vspace{8pt}

For simplicity, we assume $P$ is strongly irreducible and $\partial$-strongly irreducible and the proof for the case that $P$ is essential is the same.  By Lemma~\ref{Lessential}, either $(\partial P\cap F)\cap(\partial Q\cap F)=\emptyset$, or one can obtain an essential surface $P'$ (with $\partial P'\cap F\ne\emptyset$) by compressing and $\partial$-compressing $P$ on the same side, or after isotopy $P\cap Q$ is essential in both $P$ and $Q$ and $|\partial P\cap \partial Q|$ is minimal up to isotopy on $\partial P$ and $\partial Q$ in $\partial N$.  

If $(\partial P\cap F)\cap(\partial Q\cap F)=\emptyset$, the lemma holds trivially.  So by Lemma~\ref{Lessential}, we have the following 2 subcases to consider.

\vspace{5pt}
\noindent
\emph{Subcase 1}.  $P\cap Q$ is essential in both $P$ and $Q$, and $|\partial P\cap \partial Q|$ is minimal up to isotopy on $\partial P$ and $\partial Q$ in $\partial N$.  
\vspace{5pt}

Let $\omega=\min\{|\alpha\cap\beta|:$ $\alpha$ is a component of $\partial P\cap F$ and $\beta$ is a component of $\partial Q\cap F\}$.  Let $b_1$ and $b_2$ be the numbers of components in $P\cap F$ and $Q\cap F$ respectively.  Thus the number of arcs of $P\cap Q$ with an endpoint in $F$ is at least $\omega b_1b_2/2$.

Since $N$ contains no essential annulus with a boundary component in $F$ and since $\partial N$ is incompressible, $\chi(P)<0$ and $\chi(Q)<0$. 
As in the argument in Claim 1 above, the maximal numbers of pairwise nonparallel arcs in $P$ and $Q$ are at most $6g+3(b+b_1)-6$ and $6g+3(b+b_2)-6$ respectively.  Thus if $\omega$ is sufficiently large, there are a pair of arcs $\alpha$ and $\beta$ in $P\cap Q$ such that $\partial\alpha\cap F\ne\emptyset$, $\partial\beta\cap F\ne\emptyset$, and $\alpha$ and $\beta$ are parallel in both $P$ and $Q$.  Note that the bound for $\omega$ is an explicit quadratic function of $g$ and $b$.  Let $R_P$ and $R_Q$ be the rectangles bounded by $\alpha\cup\beta$ in $P$ and $Q$ respectively.  We can choose $\alpha$ and $\beta$ so that $\Int(R_P)\cap\Int(R_Q)=\emptyset$.  So $A=R_P\cup R_Q$ is an embedded annulus or M\"{o}bius band in $N$.  As $\alpha\cap F\ne\emptyset$, at least one component of $\partial A$ lies in $F$.  Similar to Claim 1 in Case A, since $N$ is irreducible, if $A$ is a M\"{o}bius band, $\partial A$ must be essential in $\partial N$.  Since $|\partial P\cap\partial Q|$ is minimal in their isotopy classes, if $A$ is an annulus, $\partial A$ is essential in $\partial N$.  Hence $A$ is incompressible. Since $Q$ is an essential surface and $\alpha$ is an essential arc in $Q$, $\alpha$ must be an essential arc in $N$. Hence $A$ is $\partial$-incompressible in $N$. So $A$ is essential in $N$. If $A$ is a M\"{o}bius band, a double cover of $A$ is an essential annulus.  This contradicts our hypothesis in Case B that no such essential annulus exists.  Therefore, $\omega$ and  $d(\partial P\cap F,\partial Q\cap F)$ must be bounded by a number $K'$ that depends only on $g$ and $b$.  As above, $K'$ can be chosen to be an explicit quadratic function of $g$ and $b$ and the lemma holds.

\vspace{5pt}
\noindent
\emph{Subcase 2}.  One can obtain an essential surface $P'$ (with $\partial P'\cap F\ne\emptyset$) by compressing and $\partial$-compressing $P$.
\vspace{5pt}

Since both $P'$ and $Q$ are essential, after isotopy, $P'\cap Q$ is essential in both $P'$ and $Q$ and $|\partial P'\cap \partial Q|$ is minimal up to isotopy on $\partial P'$ and $\partial Q$ in $\partial N$.  Then the argument in Subcase 1 above implies that $d(\partial P'\cap F,\partial Q\cap F)\le K'$, where $K'$ can be chosen to be an explicit quadratic function of $g$ and $b$.

By Lemma~\ref{Lcompress}, $d(\partial P'\cap F,\partial P\cap F)\le\max\{1,\ 4g+2b-2\}$. 
Thus as in the proof of Claim 2 above,  $d(\partial P\cap F,\partial Q\cap F)\le d(\partial P\cap F,\partial P'\cap F)+1+d(\partial P'\cap F,\partial Q\cap F)\le K'+\max\{1,\ 4g+2b-2\}+1$. 
\end{proof}

\begin{lemma}\label{Lbundle}
Suppose $N$ is a twisted $I$--bundle over a closed non-orientable surface and $F=\partial N$.  Let $P$ be a properly embedded orientable genus-$g$ surface with boundary  and suppose $P$ is either essential or strongly irreducible and $\partial$-strongly irreducible.   Then there is a number $K$ depending only on $g$, such that $d(\partial P,\mathcal{A}_N(F))\le K$, where $\mathcal{A}_N(F)$ is the annulus complex defined in Definition~\ref{Dann}.  Moreover, $K$ can be chosen to be an explicit linear function of $g$.
\end{lemma}
\begin{proof} Any orientable essential surface with boundary in the $I$-bundle $N$ is an annulus.  If $P$ is an essential surface, $P$ must be an annulus and $d(\partial P,\mathcal{A}_N(F))=0$.  So we may assume that $P$ is strongly irreducible and $\partial$-strongly irreducible.  

By Lemma~\ref{Lessential}, for any vertical annulus $Q$, either $\partial P\cap\partial Q=\emptyset$, or one can obtain an essential annulus by compressing and $\partial$-compressing $P$ on the same side, or after isotopy $P\cap Q$ is essential in both $P$ and $Q$.  

Suppose we can obtain an essential annulus $P'$ after some compressions and $\partial$-compressions on the same side of $P$.  By Lemma~\ref{Lcompress} and since $F=\partial N$, $P\in C_{g,0}$ and $d(\partial P,\mathcal{A}_N(F))\le d(\partial P,\partial P')\le\max\{1,\ 4g-2\}$.   Thus by Lemma~\ref{Lessential}, it remains to consider the case that for any vertical annulus $Q$ of $N$, we can isotope $P$ so that $P\cap Q$ consists of arcs essential in both $P$ and $Q$.

We may choose $Q$ to be a vertical annulus or M\"{o}bius band so that 
$\partial P\cap\partial Q$ is minimal among all vertical annuli and M\"{o}bius bands with the property that $P\cap Q$ consists of arcs essential in both $P$ and $Q$.  If $\partial P\cap\partial Q=\emptyset$ then $d(\partial P,\mathcal{A}_N(F))\le 1$.  So we may assume $\partial P\cap\partial Q\ne\emptyset$.

Let  $\omega=\min\{|\alpha\cap\beta|:$ $\alpha$ is a component of $\partial P$ and $\beta$ is a component of $\partial Q$\}.  As in the proof of Lemma~\ref{Ldist}, if $\omega$ is large, 
then there are a pair of arcs $\alpha$ and $\beta$ in $P\cap Q$ that are parallel in $P$ and we can construct a new essential annulus or M\"{o}bius band with fewer intersection with $P$.  As in the proof of Claim 1 in Lemma~\ref{Ldist}, this implies that $\omega\le\max\{5,\ 12g-6\}$.
As before, by \cite[Lemma 2.1]{MM1} and the argument at the beginning of the proof of Lemma~\ref{Ldist}, $d(\partial P,\mathcal{A}_N(F))\le d(\partial P,\partial Q)\le 2\omega +1\le \max\{11,\ 24g-11\}$ if $g(F)\ge 2$, and $d(\partial P,\mathcal{A}_N(F))\le d(\partial P,\partial Q)\le \omega \le\max\{5,\ 12g-6\}$ if $F$ is a torus.
\end{proof}

In the lemmas above, we proved some nice properties of  strongly irreducible and $\partial$-strongly irreducible surfaces in $N$.  Next we consider surfaces that are strongly irreducible but not $\partial$-strongly irreducible.

\begin{lemma}\label{Lbstr}
Let $P$ be a strongly irreducible surface properly embedded in $N$.  If $P$ is not $\partial$-strongly irreducible, then there is a surface $P'$ obtained by $\partial$-compressing $P$ and deleting any resulting $\partial$-parallel components, such that $P'$ is either 
\begin{enumerate}
\item strongly irreducible and $\partial$-strongly irreducible or 
\item essential in $N$, or
\item $P'=\emptyset$, i.e., after some $\partial$-compressions on $P$, every component of the resulting surface is $\partial$-parallel.
\end{enumerate}
\end{lemma}
\begin{proof}
Let $D$ be a $\partial$-compressing disk.  We say $D$ is disk-busting if every compressing disk on the other side of $P$ intersects $\partial D$.  If $P$ has a $\partial$-compressing disk $D$ that is not disk-busting, then we perform a $\partial$-compression along $D$.  As $D$ is not disk-busting, there is a compressing disk $D'$ of $P$ on the other side which remains a compressing disk after the $\partial$-compression.  Moreover, since $P$ is compressible on both sides, the surface obtained by $\partial$-compression along $D$ remains compressible on both sides and strongly irreducible.

After some $\partial$-compressions as above, we may assume every $\partial$-compressing disk of $P$ is disk-busting.  If $P$ is still not $\partial$-strongly irreducible, there must be a pair of $\partial$-compressing disks $D_1$ and $D_2$ on different sides of $P$ with $\partial D_1\cap\partial D_2=\emptyset$.  Now we perform $\partial$-compression on $P$ along $D_1$ and $D_2$ simultaneously and obtain a surface $P'$.  Since both $D_1$ and $D_2$ are disk-busting and $D_1$ and $D_2$ are on different sides of $P$, $P'$ is incompressible in $N$.  Therefore, after some more $\partial$-compressions on $P'$, we obtain a surface of which every component is either  essential or $\partial$-parallel.
\end{proof}

In the proof of Lemma~\ref{Lbstr}, if a component of $P$ is $\partial$-parallel and outermost, then we can simply eliminate this component.
Next we discuss how the boundary curves of $P$ change during the $\partial$-compressions in the proof of Lemma~\ref{Lbstr}. This discussion will be used later. Since we are mainly interested in the curves in $F$, we suppose all the $\partial$-compressions occur at $F$.  Each step in the proof of Lemma~\ref{Lbstr} is either a single $\partial$-compression or two simultaneous $\partial$-compressions on different sides of $P$.  The resulting surface after each step is either strongly irreducible or incompressible.  Note that a $\partial$-compression does not create any $\partial$-parallel disk, so the boundary of any resulting $\partial$-parallel component is essential in $\partial N$.  As we pointed out in Definition~\ref{D33}, the boundary curves of a strongly irreducible surface are essential in $\partial N$, so we can view them as vertices in the curve complex $\mathcal{C}(F)$.  Next we study how the distance of the boundary curves change after each step in the proof of Lemma~\ref{Lbstr}.

Let $D$ be a $\partial$-compressing disk for $P$ in the proof of Lemma~\ref{Lbstr}.  Then a $\partial$-compression along $D$ can be viewed as an isotopy pushing $D$ into a product neighborhood of $F$.  In fact, we can find a product neighborhood $F\times I$ of $F$ in $N$ such that every level $F\times\{t\}$ is transverse to $P$ except for a singular level $s\in (0,1)$ where $F\times\{s\}$ is transverse to $P$ except for a single saddle tangency.  The saddle tangency corresponds to the $\partial$-compressing disk $D$. 
Suppose $F=F\times\{0\}$ and $N'=N-(F\times[0,1))$.
Then $N'\cong N$ and $P\cap N'$ can be viewed as the surface obtained by $\partial$-compressing $P$ along $D$. 
So $\chi(P\cap(F\times I))=-1$. 
Since $P$ is orientable, the component of $P\cap(F\times I)$ that contains the saddle tangency must be a pair of pants. Hence $P\cap (F\times I)$ consists of a pair of pants and a collection of vertical annuli.  In the proof of Lemma~\ref{Lbstr}, the surface after the $\partial$-compression along $D$ remains either incompressible or strongly irreducible, so every component of $P\cap(F\times\partial I)$ is an essential curve in $F\times\partial I$.  To simplify notation, we do not distinguish a curve $\gamma$ in $F\times\{t\}$ from the vertex in $\mathcal{C}(F)$ represented by $\pi(\gamma)$, where $\pi:F\times I\to F$ is the projection.  It is  easy to see that for any curves $\gamma_0$ and $\gamma_1$ in $P\cap(F\times\{0\})$ and $P\cap(F\times\{1\})$ respectively, $d(\gamma_0,\gamma_1)\le 1=-\chi(P\cap(F\times I))$ in $\mathcal{C}(F)$. 

The situation is slightly more complicated when we simultaneously $\partial$-compressing $P$ (on different sides) along two disjoint $\partial$-compressing disks $D_1$ and $D_2$ in the last part of the proof of Lemma~\ref{Lbstr}.   
Similar to the argument above, we can find a product neighborhood $F\times I$ of $F$ in $N$ with $F\times\{0\}=F$ such that every level $F\times\{t\}$ is transverse to $P$ except for a singular level $s\in(0,1)$ where $F\times\{s\}$ is transverse to $P$ except for two saddle tangencies.  The two saddle tangencies correspond to the $\partial$-compressing disks $D_1$ and $D_2$.  Similar to the first case above, in the proof of Lemma~\ref{Lbstr}, every component of $P\cap(F\times\partial I)$ is an essential curve in $F\times\partial I$. Let $\Theta$ be the possibly disconnected graph of $P\cap (F\times\{s\})$ containing the two saddle tangencies.  So $\Theta$ has two vertices of valence 4.  If the genus of $F$ is at least 2, then there must be an essential simple closed curve $\alpha$ in $F\times\{s\}$ that is disjoint from $\Theta$ and $P\cap (F\times\{s\})$.  This implies that if $F$ is not a torus, for any components $\gamma_0$ and $\gamma_1$ of $P\cap(F\times\{0\})$ and $P\cap(F\times\{1\})$ respectively, $d(\gamma_0,\gamma_1)\le d(\gamma_0,\alpha)+d(\alpha,\gamma_1)\le 2=-\chi(P\cap(F\times I))$ in $\mathcal{C}(F)$.  

If $F$ is a torus, since $P\cap(F\times\partial I)$ consists of essential curves, each $P\cap (F\times\{i\})$ ($i=0,1$) consists of parallel curves in the torus $F$.  If $F\times\{s\}-\Theta$ is not a collection of disks, then there is an essential simple closed curve in $F\times\{s\}$ disjoint from $P$.  This implies that any curves $\gamma_0$ and $\gamma_1$ of $P\cap(F\times\{0\})$ and $P\cap(F\times\{1\})$ represent the same vertex in $\mathcal{C}(F)$ and $d(\gamma_0,\gamma_1)=0$.  Next we suppose every component of $F\times\{s\}-\Theta$ is a disk. Then $P\cap(F\times I)$ contains no vertical annulus and $P\cap(F\times I)$ can be viewed as a small neighborhood of $\Theta$. Since $P$ is separating, $P\cap(F\times\{0\})$ contains at least two curves and the argument above implies that $P\cap(F\times\{0\})$ contains exactly two curves which cut the torus $F\times\{0\}$ into two annuli $A_1$ and $A_2$.  Moreover, the two arcs $D_1\cap(F\times\{0\})$ and $D_2\cap(F\times\{0\})$ from the $\partial$-compressing disks are essential arcs in the two annuli $A_1$ and $A_2$ respectively.  So it is easy to see that $P\cap(F\times\{1\})$ also consists of exactly two curves, and for any $\gamma_0$ and $\gamma_1$ of $P\cap(F\times\{0\})$ and $P\cap(F\times\{1\})$ respectively, the intersection number of $\gamma_0$ and $\gamma_1$ (after projecting to the torus $F$) is one and hence $d(\gamma_0,\gamma_1)=1<2=-\chi(P\cap(F\times I))$ in the case that $F$ is a torus.  

Therefore, in any case, for any curves $\gamma_0$ and $\gamma_1$ of $P\cap(F\times\{0\})$ and $P\cap(F\times\{1\})$ respectively, $d(\gamma_0,\gamma_1)\le-\chi(P\cap(F\times I))$ and $-\chi(P\cap(F\times I))$ equals to the number of saddle tangencies in $F\times I$.

\section{Case I: The amalgamation surface $F$ is incompressible}\label{SI}
Let $M_1$, $M_2$, $F$ and $M=M_1\cup_F M_2$ be as in Theorem~\ref{Tmain}.  We regard $M_1$ and $M_2$ as submanifolds of $M$ with $F=\partial M_1=\partial M_2$.  In this section, we prove Theorem~\ref{Tmain} in the case that both $M_1$ and $M_2$ have incompressible boundary, i.e., $F$ is incompressible in $M$.

\begin{lemma}\label{Lesscross}
Let $M_1$, $M_2$, $M$ and $F$ be as above.  Then for any integer $g$, there is a number $K_g$ which depends only on $M_1$, $M_2$ and $g$, such that, if $d(M)>K_g$, then any closed incompressible orientable surface of genus $g$ in $M$ can be isotoped disjoint from $F$.
\end{lemma}
\begin{proof}
Let $S$ be a closed incompressible orientable surface of genus $g$ in $M$.  Suppose $S$ cannot be isotoped disjoint from $F$.  As both $S$ and $F$ are incompressible, we may assume $F\cap S$ is essential in both $F$ and $S$. 

Let $S_i=M_i\cap S$ ($i=1,2$).  So $S_i$ has no disk component, each $S_i$ is incompressible in $M_i$, $S=S_1\cup S_2$ and $\chi(S)=\chi(S_1)+\chi(S_2)$. Since $S$ cannot be isotoped disjoint from $F$, we  obtain an essential surface $S_i'$ in $M_i$ after at most $-\chi(S_i)$ $\partial$-compressions on $S_i$. Each $\partial$-compression changes the boundary curves of the surface by at most distance one in $\mathcal{C}(F)$.  Thus for any components $\gamma_i$ and $\gamma_i'$ of $\partial S_i$ and $\partial S_i'$ respectively, the distance $d(\gamma_i,\gamma_i')\le-\chi(S_i)$. 

Suppose neither $M_1$ nor $M_2$ is a twisted $I$--bundle.  Let $\Omega_i$ be the fixed essential surface with maximal Euler characteristic used in defining $d(M)$, i.e., $d(M)=d(\partial \Omega_1\cap F,\partial \Omega_2\cap F)$.  By Lemma~\ref{Ldist}, for any essential surface $Q$ properly embedded in $M_i$ with genus at most $g$ and $\partial Q\subset F$, there is a number $K_i$ such that $d(\partial \Omega_i\cap F, \partial Q)\le K_i$.  Thus there is a component $\gamma_i'$ of $\partial S_i'$, such that $d(\partial \Omega_i\cap F, \gamma_i')\le K_i$, $i=1,2$. 
Let $\gamma$ be a component of $\partial S_1=\partial S_2$.  So we have $d(\partial \Omega_1\cap F,\partial \Omega_2\cap F)\le d(\partial \Omega_1\cap F,\gamma_1')+d(\gamma_1',\gamma)+d(\gamma,\gamma_2')+d(\gamma_2',\partial \Omega_2\cap F)\le K_1-\chi(S_1)-\chi(S_2)+K_2=K_1-\chi(S)+K_2=K_1+K_2+2g-2$.
Thus Lemma~\ref{Lesscross} holds in the case that neither $M_1$ nor $M_2$ is a twisted $I$--bundle.

If $M_i$ is a twisted $I$--bundle over a closed non-orientable surface, then $S_i'$ must be a vertical annulus and each component of $\partial S_i'$ represents a vertex in the annulus complex of $M_i$.  By the definition of $d(M)$ in the case that $M_i$ is a twisted $I$--bundle, the argument above plus Lemma~\ref{Lbundle} also prove Lemma~\ref{Lesscross} in the case that some $M_i$ is a twisted $I$--bundle.
\end{proof}

Let $S$ be an unstabilized Heegaard surface of genus $g$.  As in Theorem~\ref{TST}, the untelescoping of the Heegaard splitting \cite{ST} gives a decomposition $M=N_0\cup_{F_1}N_1\cup_{F_2}\dots\cup_{F_m}N_m$, where each $F_i$ is incompressible in $M$ and $g(F_i)\le g$.  By Lemma~\ref{Lesscross}, we may assume $d(M)$ is so large that $F_i\cap F=\emptyset$ for each $i$ after isotopy.  So we may suppose $F\subset\Int(N_i)$ for some $i$.  Without loss of generality, we may assume $N_i$ is connected.  By the untelescoping construction, $N_i$ has a strongly irreducible Heegaard surface $P_i$ and $g(P_i)\le g$.  Note that if $S$ is strongly irreducible, then $N_i=M$.   The following Lemma of Bachman-Schleimer-Sedgwick \cite[Lemma 3.3]{BSS} says that we can isotope $P_i$ so that $P_i$ intersects $F$ nicely. 
If $F$ is parallel to some $F_j$ above, then Theorem~\ref{Tmain} holds.  Suppose Theorem~\ref{Tmain} is not true, then $F$ is not parallel to a component of $\partial N_i$.

\begin{lemma}[Bachman-Schleimer-Sedgwick \cite{BSS}]\label{LBSS}
Let $N_i$ be a compact, irreducible, orientable 3--manifold with $\partial N_i$ incompressible, if non-empty. Suppose $P_i$ is a strongly irreducible Heegaard surface of $N_i$.  Suppose further
that $N_i$ contains an incompressible, orientable, closed, non-boundary
parallel surface $F$. Then either
\begin{enumerate}
\item $P_i$ may be isotoped to be transverse
to $F$, with every component of $P_i - N(F)$ incompressible in the
respective submanifold of $N_i - N(F)$, where $N(F)$ is a small neighborhood of $F$ in $N_i$,
\item $P_i$ may be isotoped to be transverse
to $F$, with every component of $P_i - N(F)$ incompressible in the
respective submanifold of $N_i - N(F)$ except for exactly one
strongly irreducible component, or
\item $P_i$ may be isotoped to be almost transverse
to $F$ (i.e., $P_i$ is transverse to $F$ except for one saddle point), with every component of $P_i - N(F)$ incompressible in the
respective submanifold of $N_i - N(F)$.
\end{enumerate}
\end{lemma}

Let $N(F)=F\times I$ be a product neighborhood of $F$ in $N_i$ and let $X$ and $Y$ be the two components of $N_i-\Int(N(F))$.  As $P_i$ is a Heegaard surface of $N_i$ and the incompressible surface $F$ is not parallel to $\partial N_i$, $F\cap P_i\ne\emptyset$. Let $S_X=P_i\cap X$ and $S_Y=P_i\cap Y$. 
By Lemma~\ref{LBSS}, we may assume that each component of $S_X$ and $S_Y$ is either incompressible or strongly irreducible in $X$ and $Y$ respectively.   Moreover, both $S_X$ and $S_Y$ are essential subsurfaces of $P_i$ (i.e. $\partial S_X$ and $\partial S_Y$ are essential curves in $P_i$).  Hence $\chi(S_X)+\chi(S_Y)\ge\chi(P_i)$. By projecting $F\times I$ to $F$, we may view $\partial S_X$ and $\partial S_Y$ as curves in $F$. By Lemma~\ref{LBSS}, $P_i$ is transverse to every level surface $F\times\{t\}$ in $F\times I$ except for at most one saddle tangency which only occurs in case (3) of Lemma~\ref{LBSS}.   Thus, for any components $\gamma_X$ and $\gamma_Y$ of $\partial S_X$ and $\partial S_Y$ respectively, $d(\gamma_X, \gamma_Y)\le 1$ in $\mathcal{C}(F)$. 

Since $F\cap P_i\ne\emptyset$ after any isotopy, $S_X$ or $S_Y$ cannot be changed to a set of $\partial$-parallel surfaces by $\partial$-compressions on $S_X$ or $S_Y$ in $X$ or $Y$ respectively.  Thus by Lemma~\ref{Lbstr}, we can obtain a pair of surfaces $S_X'$ and $S_Y'$ by some $\partial$-compressions on $S_X$ and $S_Y$ respectively, such that $S_X'$ and $S_Y'$ are either essential or strongly irreducible and $\partial$-strongly irreducible in $X$ and $Y$ respectively.  The numbers of $\partial$-compressions on $S_X$ and $S_Y$ are at most $-\chi(S_X)$ and $-\chi(S_Y)$ respectively.  Since each $\partial$-compression changes a curve by distance at most one in the $\mathcal{C}(F)$, by the argument after Lemma~\ref{Lbstr}, for any components $\gamma_X$ and $\gamma_X'$ of $\partial S_X$ and $\partial S_X'$ respectively, $d(\gamma_X,\gamma_X')\le-\chi(S_X)$.  

Since $-\chi(S_X)-\chi(S_Y)\le -\chi(P_i)\le 2g-2$ and since $d(\gamma_X, \gamma_Y)\le 1$ for any components $\gamma_X$ and $\gamma_Y$ of $\partial S_X$ and $\partial S_Y$ respectively, for any components $\gamma_X'$ and $\gamma_Y'$ in $\partial S_X'$ and $\partial S_Y'$, we have $d(\gamma_X', \gamma_Y')\le d(\gamma_X',\gamma_X)+d(\gamma_X,\gamma_Y)+d(\gamma_Y,\gamma_Y')\le -\chi(S_X)+1-\chi(S_Y)\le 1-\chi(P_i)\le 2g-1$.

Note that since $N_i$ is a submanifold of $M=M_1\cup_F M_2$ with $F\subset\Int(N_i)$, to simplify notation, we will regard $X$ and $Y$ as submanifolds of $M_1$ and $M_2$ respectively with $F\subset \partial X$  and $F\subset \partial Y$.   Since $F$ is not parallel to a component of $\partial N_i$, $X$ and $Y$ are not $I$--bundles unless $M_1$ or $M_2$ is a twisted $I$--bundle.

We first suppose neither $M_1$ nor $M_2$ is a twisted $I$--bundle. 
Let $\Omega_j$ ($j=1,2$) be the fixed essential surface in $M_j$ used in defining $d(M)$, i.e., $d(M)=d(\partial\Omega_1\cap F,\partial\Omega_2\cap F)$.
Since $\partial X$ and $\partial Y$ are incompressible in $M_1$ and $M_2$ respectively, we may assume $\Omega_1\cap X$ and $\Omega_2\cap Y$ are essential surfaces in $X$ and $Y$ respectively.  Moreover, we may assume $\Omega_1\cap X$ and $\Omega_2\cap Y$ are essential subsurfaces of $\Omega_1$ and $\Omega_2$ respectively, and in particular, 
$\chi(\Omega_1\cap X)\ge\chi(\Omega_1)$ and $\chi(\Omega_2\cap Y)\ge\chi(\Omega_2)$.  

As $F\cap \partial (\Omega_1\cap X)=\partial\Omega_1\cap F$ and $F\cap \partial (\Omega_2\cap Y)=\partial\Omega_2\cap F$, 
by applying Lemma~\ref{Ldist} to $X$ and $Y$, we conclude that there is a number $K$ depending only on $g$ and $\max\{-\chi(\Omega_1),-\chi(\Omega_2)\}$, such that $d(\partial S_X',\partial\Omega_1\cap F)\le K$ and $d(\partial S_Y',\partial\Omega_2\cap F)\le K$.  Let $\gamma_X'$ and $\gamma_Y'$ be any components of $\partial S_X'$ and $\partial S_Y'$ respectively.  Recall that we have concluded earlier that $d(\gamma_X',\gamma_Y')\le 2g-1$.  So we have  
$d(M)=d(\partial\Omega_1\cap F,\partial\Omega_2\cap F)\le d(\partial\Omega_1\cap F,\gamma_X') + d(\gamma_X',\gamma_Y')+ d(\gamma_Y',\partial\Omega_2\cap F)\le K+(2g-1)+K=2K+2g-1$.

If $M_j$ is a twisted $I$--bundle, then it is possible that $X=M_1$ or $Y=M_2$ is a twisted $I$--bundle. In this case, we can replace $\partial\Omega_j\cap F$ by the annulus complex $\mathcal{A}_F(M_j)$ in the argument above.  We can apply Lemma~\ref{Lbundle} instead of Lemma~\ref{Ldist} and get the same inequalities on $d(M)$.  

Therefore, if $d(M)$ is sufficiently large, we get a contradiction and this means that $F$ must be parallel to some incompressible surface $F_i$ in the untelescoping and Theorem~\ref{Tmain} holds in the case that both $M_1$ and $M_2$ have incompressible boundary.

\section{Case II: The amalgamation surface $F$ is compressible on both sides}\label{SII}

The case that both $M_1$ and $M_2$ have compressible boundary in Theorem~\ref{Tmain} basically follows from a theorem of Scharlemann and Tomova \cite{ST1} and  a theorem
of Hartshorn \cite{H}, also see \cite{L3}.

Let $\mathcal{D}_i$ be the disk complex of $M_i$ ($i=1,2$).  Recall that in this case $d(M)$ is defined to be $d(\mathcal{D}_1, \mathcal{D}_2)$.  We may assume $d(\mathcal{D}_1, \mathcal{D}_2)\ge 2$ which implies that $F$ is strongly irreducible in $M$.  By Casson-Gordon \cite{CG} and Haken's lemma \cite{Hak}, this also implies that $M=M_1\cup_F M_2$ is irreducible and $\partial$-irreducible and $M$ is not $S^3$.  Hence Theorem~\ref{TS3} holds in this case.  

Next we suppose $d(M)>2g$.   Let $S$ be an unstabilized Heegaard surface of genus $g$.  As in Theorem~\ref{TST}, the untelescoping of the Heegaard splitting \cite{ST} gives a decomposition $M=N_0\cup_{F_1}N_1\cup_{F_2}\dots\cup_{F_m}N_m$, where each $F_i$ is incompressible in $M$ and $g(F_i)\le g$.  By Hartshorn's theorem \cite{H}, see \cite{L3} for another proof, either $F_i\cap F=\emptyset$ after isotopy or $d(M)=d(\mathcal{D}_1, \mathcal{D}_2)\le 2g(F_i)\le 2g$ for each $i$.  Since $d(M)>2g$, we may assume $F\subset\Int(N_k)$ for some $k$.  Without loss of generality, we suppose $N_k$ is connected.  As in Theorem~\ref{TST}, there is a strongly irreducible Heegaard surface $P_k$ of the 3--manifold $N_k$ and $g(P_k)\le g$.

Let $Q_j$ ($j=1,2$) be the surface obtained by maximally compressing $F$ in $M_j$ and removing all resulting 2--sphere components.  We may assume $Q_j\subset\Int(M_j)$, $Q_1\cup Q_2$ bounds a submanifold $M_F$ in $M$, and $F$ is a strongly irreducible Heegaard surface of $M_F$.  Since $F\subset\Int(N_k)$ and $\partial N_k$ is incompressible in $M$, any compressing disk for $F$ can be isotoped into $N_k$.  So after isotopy, we may assume $M_F\subset N_k$.

Since $P_k$ is strongly irreducible, a theorem of Scharlemann and Tomova \cite{ST1} says that either $d(\mathcal{D}_1, \mathcal{D}_2)\le 2g(P_k)\le 2g$, or $F$ and $P_k$ are well-separated, or $F$ and $P_k$ are parallel.  

Next we show that $F$ and $P_k$ are not well-separated.  Suppose on the contrary that they are well-separated, i.e., $M_F$ can be isotoped disjoint from $N_k$.  Let $M_F'$ be a submanifold of $M$ that is isotopic to $M_F$ and disjoint from $N_k$.  Since $M_F\subset N_k$ and $M_F'\cap N_k=\emptyset$, $M_F'$ is disjoint from $M_F$.  Recall that $M=M_1\cup_F M_2$ is an amalgamation of $M_1$ and $M_2$ along $F$, so $M_F'$ lies in either $\Int(M_1)$ or $\Int(M_2)$. Without loss of generality, suppose $M_F'\subset\Int(M_1)$.  By our construction of $M_F$, the surface $Q_2\subset\partial M_F$ lies in $M_2$.  Let $Y$ be a component of $Q_2$ and let $Y'$ be the component of $\partial M_F'$ isotopic to $Y$.  $Y$ and $Y'$ are two-sided, incompressible, disjoint and isotopic surfaces in $M$, so $Y\cup Y'$ bounds a product region $Y\times I$ in $M$ (one can see this easily after lifting $Y$ and $Y'$ to the covering space of $M$ corresponding to $\pi_1(Y)$).  Since $Y\subset \Int(M_2)$ and $Y'\subset\partial M_F'\subset\Int(M_1)$, the product region $Y\times I$ must contain the amalgamation surface $F$.  Moreover, since $Y$ and $Y'$ are incompressible in $M$ and $F\subset Y\times I$, a compressing disk for $F$ can be isotoped into $Y\times I$.  Hence after isotopy, we may assume $M_F$ lies in the production region $Y\times I$.  Each closed incompressible surface in the product $Y\times I$ is parallel to $Y$.  This implies that $M_F$ is isotopic to $Y\times I$ and $F$ can be viewed as a Heegaard surface of $Y\times I$.  By \cite{ST2},  Heegaard splittings of a product $Y\times I$ are standard and in particular the distance of the Heegaard splitting of $Y\times I$ along $F$ is at most $2$.  This contradicts our assumption that $d(M)>2g\ge 2$.  So $F$ and $P_k$ are not well-separated.

So the theorem of Scharlemann and Tomova \cite{ST1} implies that if $d(M)>2g$, then $F$ and $P_k$ must be parallel.  Therefore, Theorem~\ref{Tmain} holds if both $M_1$ and $M_2$ have compressible boundary, and in this case we may choose the bound $K=2g$.

\section{Case III: The amalgamation surface $F$ is compressible on one side}\label{Scase3}
In the next two sections, we suppose $F$ is compressible in $M_1$ but incompressible in $M_2$.   We denote the disk complex of $M_1$ by $\mathcal{D}_1$.

\begin{proposition}\label{Ptrivial}
Let $\gamma$ be a nontrivial simple closed curve in $F$.  Suppose $\gamma$ bounds an embedded disk in $M=M_1\cup_F M_2$.  Then $d(\gamma, \mathcal{D}_1)\le 1$.
\end{proposition}
\begin{proof}
Let $D$ be the embedded disk bounded by $\gamma$ in $M$.  We may assume that $|D\cap F|$ is minimal among all disks bounded by $\gamma$ and  transverse to $F$.  Since $F$ is incompressible in $M_2$, if $\Int(D)\cap F=\emptyset$ then $D$ must be a compressing disk of $M_1$ and $d(\gamma, \mathcal{D}_1)=0$. 

Let $\gamma'$ be a component of $D\cap F$ that is innermost in $D$ and let $\delta$ be the subdisk of $D$ bounded by $\gamma'$.  If $\gamma'$ is a trivial curve in $F$, then a standard cutting and pasting yields a new disk bounded by $\gamma$ with fewer intersection curves with $F$.  So $\delta$ must be a compressing disk in $M_1$.  Since $D$ is embedded, $\gamma$ and $\gamma'$ are disjoint in $F$.  Therefore, $d(\gamma, \mathcal{D}_1)\le d(\gamma, \gamma')\le 1$.
\end{proof}

Many parts of the proof of Theorem~\ref{Tmain} in this case comes down to the situation that $F$ lies in a submanifold $M'$ of $M$ with incompressible boundary, and we need to study how various surfaces intersect $F$.  The following technical lemma deals with this situation and will be used in several places of the proof.  The key point of Lemma~\ref{Lpartial} is that the bound on the distance in Lemma~\ref{Lpartial} depends on $-\chi(P)$ not on $-\chi(P_2)$ which can be large since $F$ is compressible in $M_1$.

\begin{lemma}\label{Lpartial}
Let $M'$ be a compact submanifold of $M=M_1\cup_F M_2$ with $F\subset\Int(M')$ and suppose $\partial M'$ is incompressible in $M'$. Let $P$ be an orientable connected surface properly embedded in $M'$.  
Suppose $P$ is either incompressible or strongly irreducible in $M'$, $P\cap F\ne\emptyset$, and each component of $P\cap F$ is essential in $F$. Let $M_2'=M_2\cap M'$ and $P_2=P\cap M_2'$. Suppose $P_2$ is either incompressible or strongly irreducible in $M_2'$ and $P_2$ does not lie in a product neighborhood of $F$ in $M_2'$. 
Then there is a surface $Q$ obtained by some $\partial$-compressions on $P_2$ in $M_2'$ and removing all resulting $\partial$-parallel components, such that $d(Q\cap F,(P\cap F)\cup\mathcal{D}_1)\le\max\{3-\chi(P), 2\}$ and $Q$ is either an essential or a strongly irreducible and $\partial$-strongly irreducible surface properly embedded in $M_2'$.
\end{lemma}
\begin{proof}
If $P_2$ is incompressible in $M_2'$, then after performing some $\partial$--compressions on $P_2$ in $M_2'$, we get a surface $Q$ such that each component of $Q$ is either essential or  $\partial$--parallel in $M_2'$.  Similarly, if $P_2$ is strongly irreducible but not $\partial$-strongly irreducible, as in Lemma~\ref{Lbstr}, we can obtain a surface $Q$ after some $\partial$-compressions on $P_2$ in $M_2'$ such that each component of $Q$ is either essential, or $\partial$-parallel,  or strongly irreducible and $\partial$-strongly irreducible in $M_2'$.   Since $P_2$ does not lie in a product neighborhood of $F$, after discarding all the $\partial$--parallel components, we get a surface $Q$ which is either essential or strongly irreducible and $\partial$-strongly irreducible in $M_2'$.  As $\partial P_2\cap F=P\cap F$, to prove the lemma, we  need to study the distance between $\partial P_2$ and $\partial Q$ in the curve complex $\mathcal{C}(F)$.

The surface $F$ is a boundary component of $M_2'$.  Since we are only interested in how the curves in $\partial P_2\cap F$ change during $\partial$-compressions, to simplify notation, we will assume that all the $\partial$-compressions on $P_2$ in the construction above occur at $F$, i.e., for any $\partial$-compressing disk $D$ for $P_2$, we assume $D\cap\partial M_2'\subset F$.

A $\partial$-compression on $P_2$ is basically the same as an isotopy that pushes the $\partial$-compressing disk into a product neighborhood of $F$.  Thus we can find a product neighborhood $F\times I$ of $F$ in $M_2'$ and assume $Q=P_2\cap\overline{M_2'-(F\times I)}$.  We denote $F\times\{t\}$ by $F_t$ and suppose $F_0=F\subset\partial M_2'$.  By the discussion after the proof of Lemma~\ref{Lbstr}, we may describe each $\partial$-compression using a saddle tangency in $F_t\cap P_2$.  In the proof of Lemma~\ref{Lbstr}, we have to simultaneously perform two $\partial$-compressions, so we allow two saddle tangencies at the same level surface $F_t$.  Since $Q$ is obtained by a sequence of $\partial$-compressions and pushing away the $\partial$-parallel components, we may assume that there are finitely many numbers $0=s_0<s_1<\cdots < s_k=1$, such that
\begin{enumerate}
\item $P_2$ is transverse to each $F_{s_i}$, and each component of $P_2\cap F_{s_i}$ is essential in $F_{s_i}$, 
\item for each $i$, there is one special component of $P_2\cap (F\times [s_i,s_{i+1}])$ that is transverse to every $F_t$ except for a singular level $t_i\in (s_i,s_{i+1})$ where it is transverse to $F_{t_i}$ except for one or two saddle tangencies, and 
\item every other component of $P_2\cap (F\times [s_i,s_{i+1}])$ is either a vertical annulus or a $\partial$-parallel surface in $F\times [s_i,s_{i+1}]$ with boundary in $F_{s_i}$. 
\end{enumerate}
Each saddle tangency in the special component in (2) corresponds to a $\partial$-compression on $P_2$ and the $\partial$-parallel components in (3) are the possible $\partial$-parallel components after a $\partial$-compression. Note that it is possible to have two saddle tangencies at the same level $F_{t_i}$ because in the proof of Lemma~\ref{Lbstr}, we have to simultaneously $\partial$-compressing the surface on both sides in order to obtain an incompressible surface, see the discussion after the proof of Lemma~\ref{Lbstr}.  We regard $Q=P_2\cap\overline{M_2'-(F\times I)}$, so $\partial Q\subset F_1\cup(\partial M_2'-F_0)$.

To simplify notation, we do not distinguish a nontrivial curve $\gamma$ in $F_t$ from the vertex in $\mathcal{C}(F)$ representing $\pi(\gamma)$, where $\pi:F\times I\to F$ is the projection.  Next we show that $d(Q\cap F_1,(P\cap F_0)\cup\mathcal{D}_1)\le\max\{3-\chi(P), 2\}$. The argument is similar to \cite[Claims 1 and 3 of Lemma 2.2]{L3}.

We first point out a useful fact on curves in $P_2\cap F_{s_i}$. 
Let $\gamma_i$ be any component of $P_2\cap F_{s_i}$ and let $Q_i$ be the component of $P_2\cap (F\times [s_{i-1},s_{i}])$ that contains $\gamma_i$.  By our assumption on $P_2\cap (F\times [s_{i-1},s_{i}])$ above, $Q_i$ is either a vertical annulus or a special component in (2) above.  Since the saddle tangencies in a special component correspond to $\partial$-compressions, $\partial Q_i\cap F_{s_{i-1}}\ne\emptyset$.  Let $\gamma_{i-1}$ be a component of $\partial Q_i\cap F_{s_{i-1}}$.  By the discussion after the proof of Lemma~\ref{Lbstr}, $d(\gamma_{i-1},\gamma_i)\le n_i$, where $n_i$ is the number of saddle tangencies in the special component of $P_2\cap (F\times [s_{i-1},s_{i}])$ and $n_i$ is either $1$ or $2$.  
Thus we can successively find a curve $\gamma_i$ in each $P_2\cap F_{s_i}$ such that $d(\gamma_{i-1},\gamma_i)\le n_i$ for each $i$, where $n_i=1$ or 2 is the number of saddle tangencies in the special component of $P_2\cap (F\times [s_{i-1},s_{i}])$.

We say a component $\gamma$ of $P_2\cap F_{s_i}$ is \emph{good} if the component of $P_2\cap(F\times[s_i,1])$, denoted by $Q_\gamma$, that contains $\gamma$ has a boundary component in $F_1$, i.e. $Q_\gamma\cap F_1\ne\emptyset$.  Moreover, every component of $P_2\cap F_1=Q\cap F_1$ is regarded as a good component.  Let $C_i$ be the set of good components of $P_2\cap F_{s_i}$. As $s_k=1$, $C_k=Q\cap F_1$.  Since $P_2$ does not lie in a product neighborhood of $F$ in $M_2'$, $C_i\ne\emptyset$ for all $i$.

Suppose the lemma is not true and $d(Q\cap F_1, (P\cap F_0)\cup\mathcal{D}_1)>2$. Since $s_k=1$ and $C_k=Q\cap F_1$, we have $d(C_k, (P\cap F_0)\cup\mathcal{D}_1)>2$. As $s_0=0$ and $P_2\cap F_{s_0}=P\cap F_0\supset C_0$, we have $d(C_0, P\cap F_0)\cup\mathcal{D}_1)=0$.
Let $m$ be the smallest number ($1\le m\le k$) such that $d(C_m, (P\cap F_0)\cup\mathcal{D}_1)\ge 2$.  Since $m$ is the smallest such number and $m\ge1$, $d(C_{m-1},(P\cap F_0)\cup\mathcal{D}_1)\le 1$.
By the discussion after the proof of Lemma~\ref{Lbstr} and as above, for any curves $\alpha$ and $\beta$ in $C_{m-1}$ and $C_m$ respectively, $d(\alpha,\beta)$ is smaller than or equal to the number of saddle tangencies in the special component of $P_2\cap(F\times[s_{m-1},s_m])$ and $d(\alpha,\beta)\le 2$.  We may choose $\alpha$ to be the curve in $C_{m-1}$ realizing $d(\alpha,(P\cap F_0)\cup\mathcal{D}_1)=d(C_{m-1},(P\cap F_0)\cup\mathcal{D}_1)$.  So $d(\alpha,(P\cap F_0)\cup\mathcal{D}_1)\le 1$.  Hence for any curve $\beta$ in $C_m$, $d(\beta, (P\cap F_0)\cup\mathcal{D}_1)\le d(\beta, \alpha)+ d(\alpha,(P\cap F_0)\cup\mathcal{D}_1)\le 2+1=3$.

Let $Q'$ be a component of $P_2\cap (F\times[s_m, 1])$ connecting $F_{s_m}$ and $F_1$, i.e. $\partial Q'$ contains curves in both $F_{s_m}$ and $F_1$. By the definition of $C_i$, $\partial Q'\cap F_{s_m}\subset C_m$ and $\partial Q'\cap F_1\subset C_k$.  Since $d(C_m, (P\cap F_0)\cup\mathcal{D}_1)\ge 2$, we have $d(C_m, \mathcal{D}_1)\ge 2$.  Similarly, $d(C_k,\mathcal{D}_1)\ge 2$ by our assumption. Hence $d(\partial Q', \mathcal{D}_1)\ge 2$.  This implies that $\partial Q'$ are essential curves in $P$, to see this, if a curve $\gamma$ in $\partial Q'$ is trivial in $P$, then $\gamma$ bounds a disk in $P$ and by Proposition~\ref{Ptrivial}, $d(\partial Q', \mathcal{D}_1)\le d(\gamma, \mathcal{D}_1)\le 1$, a contradiction.  Thus $Q'$ must be an essential subsurface of $P$, and $P$ cannot be a 2--sphere or disk.  In particular, $\chi(P)\le \chi(Q')\le 0$. 

Let $\Gamma$ be the total number of saddle tangencies in those special components of $Q'\cap (F\times[s_i,s_{i+1}])$, $i=m,\dots, k-1$.  Note that we are only counting the saddle tangencies in $Q'$ not all saddle tangencies.  As $\partial Q'\subset F_{s_m}\cup F_{s_k}$ ($s_k=1$), by our construction, $-\chi(Q')\ge\Gamma$ (note that this is an inequality because a component of $Q'\cap (F\times[s_i,s_{i+1}])$ may be $\partial$-parallel as in part (3) of our assumption above).  Hence $-\chi(P)\ge-\chi(Q')\ge\Gamma$.

Let $\gamma_k$ be a component of $\partial Q'\cap F_{s_k}\subset C_k=Q\cap F_1$.
Recall that in the argument earlier, we showed that we can successively find a curve $\gamma_i$ in each $Q'\cap F_{s_i}$ ($i=1,\dots,k$) such that $d(\gamma_{i-1},\gamma_i)\le n_i$ for all $i$, where $n_i$ is the number of saddle tangencies in the special component of $Q'\cap (F\times [s_{i-1},s_{i}])$.  Thus $d(\gamma_m,\gamma_k)\le \Gamma\le-\chi(Q')\le-\chi(P)$. Recall that we have concluded earlier that $d(\gamma_m,(P\cap F_0)\cup\mathcal{D}_1)\le 3$. Since $s_k=1$ and $\gamma_k$ is a component of $P_2\cap F_{s_k}=Q\cap F_{1}$,  we have
$$d(Q\cap F_1, (P\cap F_0)\cup\mathcal{D}_1)\le d(\gamma_k, \gamma_m)+d(\gamma_m, (P\cap F_0)\cup\mathcal{D}_1)\le \Gamma+3\le 3-\chi(P).$$
\end{proof}

\begin{corollary}\label{Cirred}
There is a number $K$ depending on $M_2$ such that if $d(M)\ge K$ then $M=M_1\cup_F M_2$ is irreducible and $\partial$-irreducible. 
\end{corollary}
\begin{proof}
Suppose $M$ is reducible or $\partial$-reducible and let $P$ be either an essential 2--sphere or a compressing disk in $M$.  Since $M_i$ is irreducible and $\partial M_i-F$ is incompressible in $M_i$ for both $i=1,2$, $P\cap F\ne\emptyset$ after any isotopy on $P$.  If $P\cap M_2$ is compressible in $M_2$, then we can compress $P\cap M_2$ and obtain a new essential 2--sphere or compressing disk in $M$.  So after finitely many such operations, we may assume $P\cap M_2$ is incompressible in $M_2$.
Moreover, as in Lemma~\ref{Lpartial}, after pushing parts of $P\cap M_2$ into $M_1$ via $\partial$-compressions, we may assume that $Q=P\cap M_2$ is an essential planar surface in $M_2$.  Since each component of $\partial Q$ bounds a disk in $P$, by Proposition~\ref{Ptrivial}, $d(\gamma, \mathcal{D}_1)\le 1$ for each component $\gamma$ of $\partial Q$.  

If $M_2$ is not a twisted $I$--bundle, let $\Omega_2$ be the fixed essential surface in $M_2$ used in defining $d(M)=d(\mathcal{D}_1, \partial\Omega_2\cap F)$.  Since $Q$ is planar with all but at most one boundary component in $F$, by Lemma~\ref{Ldist}, there is a number $K'$ depending on $\Omega_2$ such that $d(\partial\Omega_2\cap F,\gamma)\le K'$, where $\gamma$ is a component of $\partial Q$.  Hence $d(M)=d(\partial\Omega_2\cap F,\mathcal{D}_1)\le d(\partial\Omega_2\cap F,\gamma)+ d(\gamma, \mathcal{D}_1)\le K'+1$. 

If $M_2$ is a twisted $I$--bundle, then $Q$ must be an essential annulus and hence $d(M)=d(\mathcal{A}_{M_2},\mathcal{D}_1)\le 1$, where $\mathcal{A}_{M_2}$ is the annulus complex of the twisted $I$--bundle.

Thus in any case, if $d(M)>K'+1$, no essential 2-sphere or compressing disk $P$ exists. 
\end{proof}

\begin{corollary}\label{Cincomp}
Let $F'$ be the surface obtained by maximally compressing $F$ in $M_1$ and removing all resulting 2--sphere components. Suppose $F'\ne\emptyset$.  Then there is a number $K$ depending on $M_2$ such that if $d(M)\ge K$, $F'$ is incompressible in $M$.
\end{corollary}
\begin{proof} We may assume that $F'\subset\Int(M_1)$. By our construction, $F'$ is incompressible in $M_1$.  Suppose $F'$ is compressible in $M$ and let $D$ be a compressing disk.  So $D\cap F\ne\emptyset$. As in Corollary~\ref{Cirred}, we may assume a component $Q$ of $D\cap M_2$ is essential in $M_2$.  By Proposition~\ref{Ptrivial}, $d(\gamma, \mathcal{D}_1)\le 1$ for each component $\gamma$ of $\partial Q$.  Now the proof is the same as the proof of Corollary~\ref{Cirred}.
\end{proof}

Part (2) of Theorem~\ref{TS3} also follows from the arguments above.  However, we also need the following theorem from \cite{L6} which says that a graph complement in $S^3$ always contains a nice planar surface.

\begin{lemma}[\cite{L6}]\label{Lthin}
Let $\Gamma$ be any graph in $S^3$. Then there is a planar surface $P$ properly embedded in $S^3-N(\Gamma)$ such that all but at most one of the components of $\partial P$ bound compressing disks in the handlebody $\overline{N(\Gamma)}$ and $P$ is either 
\begin{enumerate}
\item  strongly irreducible and $\partial$-strongly irreducible, or
\item  essential (possibly an essential disk), or
\item  nonseparating and incompressible in $S^3-N(\Gamma)$. 
\end{enumerate}
\end{lemma}

\begin{proof}[Proof of Theorem~\ref{TS3}]  If $F$ is incompressible in both $M_1$ and $M_2$, then Theorem~\ref{TS3} holds trivially. If $F$ is compressible in both $M_1$ and $M_2$, then as in section~\ref{SII}, Casson-Gordon \cite{CG} implies that if $d(M)\ge 2$, then $M$ is irreducible and $\partial$-irreducible and $M\not\cong S^3$.  Therefore we only need to consider the case that $F$ is compressible on one side.  Suppose $F$ is compressible in $M_1$ but incompressible in $M_2$.  Part (1) of Theorem~\ref{TS3} in this case is proved in Corollary~\ref{Cirred}.  Suppose part (2) of Theorem~\ref{TS3} fails and $M=S^3$.

Since $M=S^3$ does not contain an incompressible surface, by Corollary ~\ref{Cincomp}, if $d(M)$ is sufficiently large, we may suppose $F'=\emptyset$ and $M_1$ must be a handlebody. We may view $M_1$ as a neighborhood of a graph in $M=S^3$.  So there is a planar surface $P$ properly embedded in $M_2$ as in Lemma~\ref{Lthin}.  As $F$ is incompressible in $M_2$, $P$ is not a compressing disk in $M_2$ and hence a component of $\partial P$ bounds a compressing disk in $M_1$ and $d(\mathcal{D}_1,\partial P)=0$.  

If $P$ is nonseparating and incompressible as in part (3) of Lemma~\ref{Lthin}, then one can perform some $\partial$-compressions and obtains an essential planar surface $Q$.  Moreover, by Lemma~\ref{Lcompress}, $d(\partial P,\partial Q)\le\max\{1,\ 4g+2b-2\}=1$ since $g=0$ and $b=0$ in this case.  Since $d(\mathcal{D}_1,\partial P)=0$, this means that $d(\mathcal{D}_1,\partial Q)\le 2$.  Thus in any of the 3 possibilities of Lemma~\ref{Lthin}, we have a planar surface $Q$ in $M_2$ that is either essential or strongly irreducible and $\partial$-strongly irreducible such that $d(\mathcal{D}_1,\partial Q)\le 2$.

Since $M=S^3$, $M_2$ cannot be a twisted $I$--bundle. Let $\Omega_2$ be the fixed essential surface in $M_2$ used in defining $d(M)$.  Since $Q$ is planar, by Lemma~\ref{Ldist}, there is a number $K'$ depending on $g(\Omega_2)$ such that $d(\partial\Omega_2\cap F,\gamma)\le K'$, where $\gamma$ is a component of $\partial Q$.  Hence $d(M)=d(\partial\Omega_2\cap F,\mathcal{D}_1)\le d(\partial\Omega_2\cap F,\gamma)+ d(\gamma, \mathcal{D}_1)\le K'+2$.  Thus if $d(M)>K'+2$, $M\not\cong S^3$.
\end{proof}

In the remainder of the paper, we assume $M$ is not $S^3$ and hence our Heegaard surfaces are not $S^2$.

\begin{lemma}\label{Llarge}
For any $g\ge 1$, there is a number $K$ depending only on $M_2$ and $g$, such that if $d(M)\ge K$ then any closed orientable incompressible surface in $M$ of genus $g$ can be isotoped disjoint from $F$.  
\end{lemma}
\begin{proof}
By Corollary~\ref{Cirred}, we may assume $d(M)$ is so large that $M$ is irreducible and $\partial$-irreducible.  Let $P$ be a closed orientable incompressible surface in $M$ of genus $g$ and suppose $F\cap P\ne\emptyset$ after any isotopy.  

Let $D$ be a compressing disk for $F$ in $M_1$.  If $P\cap D$ contains a closed curve, since $P$ is incompressible, a standard isotopy on $P$ can remove this intersection curve.  Moreover, by shrinking $D$ to be sufficiently small while fixing $P$, we can also isotope $F$ to eliminate all the arcs in $P\cap D$.  Thus, after isotopy, we may assume $P\cap D=\emptyset$.  Since $P$ is incompressible and $M$ is irreducible, after isotopy, we may also assume every curve in $P\cap F$ is essential in $F$.  Since $D\cap P=\emptyset$, for any component $\gamma$ of $P\cap F$, $d(\gamma,\mathcal{D}_1)\le d(\gamma,\partial D)\le 1$.

Since $P$ is incompressible in $M$ and $M$ is irreducible, after some isotopy, we may assume that $P\cap M_2$ is incompressible in $M_2$. 
Now we apply Lemma~\ref{Lpartial}, setting $M'$, $P$ and $P_2$ in Lemma~\ref{Lpartial} to be $M$, $P$ and $P\cap M_2$ above respectively. 
 By Lemma~\ref{Lpartial} and since $F\cap P\ne\emptyset$ after any isotopy, there is an essential surface $Q$ in $M_2$ obtained by $\partial$-compressing $P\cap M_2$ such that $d(\partial Q, (P\cap F)\cup\mathcal{D}_1)\le 3-\chi(P)=2g+1$.  For any component $\gamma$  of $P\cap F$, by our earlier assumption, $d(\gamma,\mathcal{D}_1)\le 1$.  This implies that $d(\partial Q, \mathcal{D}_1)\le 2g+2$.  Moreover, by our construction, the genus of $Q$ is at most $g$.  

If $M_2$ is a twisted $I$--bundle, then $Q$ must be an essential annulus and hence $d(M)=d(\mathcal{A}_{M_2},\mathcal{D}_1)\le d(\partial Q, \mathcal{D}_1)\le 2g+2$, where $\mathcal{A}_{M_2}$ is the annulus complex of the twisted $I$--bundle.

If $M_2$ is not a twisted $I$--bundle, let $\Omega_2$ be the fixed essential surface in $M_2$ used in defining $d(M)=d(\mathcal{D}_1, \partial\Omega_2\cap F)$.  Since $g(Q)\le g$, by Lemma~\ref{Ldist}, there is a number $K'$ depending on $\Omega_2$ and $g$, such that $d(\partial\Omega_2\cap F,\partial Q)\le K'$.  Hence we can find a  component $\gamma_Q$ of $\partial Q$ such that $d(M)=d(\partial\Omega_2\cap F,\mathcal{D}_1)\le d(\partial\Omega_2\cap F,\gamma_Q)+d(\gamma_Q,\mathcal{D}_1)\le K'+2g+2$. 
Thus if $d(M)>K'+2g+2$, $P\cap F=\emptyset$ after isotopy.
\end{proof}

Let $S$ be an unstabilized Heegaard surface of genus $g$.  As in Theorem~\ref{TST}, the untelescoping of the Heegaard splitting \cite{ST} gives a decomposition $M=N_0\cup_{F_1}N_1\cup_{F_2}\dots\cup_{F_m}N_m$, where each $F_i$ is incompressible in $M$ and $g(F_i)\le g$. By Lemma~\ref{Llarge}, we may assume $d(M)$ is so large that each $F_i$ is disjoint from $F$ after some isotopy.  Thus we may assume $F\subset N_j$ for some $j$.  Without loss of generality, we may suppose $N_j$ is connected.  Now we consider the strongly irreducible Heegaard surface $P_j$ of $N_j$ in the untelescoping construction.  Let $X$ and $Y$ be the two compression bodies in the Heegaard splitting of $N_j$ along $P_j$.  We have $P_j=\partial_+X=\partial_+Y$ and $g(P_j)\le g$. 

Let $F'$ be the surface obtained by maximally compressing $F$ in $M_1$ and removing all resulting 2--sphere components.  We may assume $F'\subset\Int(M_1)$. Let $M_F$ be the compression body bounded by $F$ and $F'$ in $M_1$.  If $F'=\emptyset$ then $M_F=M_1$ is a handlebody.

In the next two lemmas, we prove that if $d(M)$ is sufficiently large, then $F$ cannot lie in a product neighborhood of any incompressible surface $F_i$ in the untelescoping. 

\begin{lemma}\label{Lknotted}
Let $E$ be an orientable incompressible closed surface in $M$ and $E\times I$ a product neighborhood of $E$.  Suppose $M_2\subset\Int(E\times I)$, then $d(M)<K$ for some $K$ depending only on $M_2$.
\end{lemma}
\begin{proof} The hypothesis that $M_2\subset\Int(E\times I)\subset\Int(M)$ implies that $F\subset\Int(E\times I)$ and $\partial M_2=F$ (i.e. $M_2$ has no other boundary component).   
Since $E$ is incompressible in $M$ and $F\subset\Int(E\times I)$, every compressing disk for $F$ can be isotoped into $E\times I$.  Thus, after isotopy, we may assume the compression body $M_F$ described above lies in $E\times I$.  As $M_2\subset E\times I$ and $E\times I\ne M$, $F'\ne\emptyset$.  

By our construction of $M_F$ and since $\partial M_2=F$, the submanifold $M_2\cup M_F$ of $M$ is bounded by $F'$.  The assumption above says that $M_2\cup M_F\subset E\times I$.  
By Corollary~\ref{Cincomp}, we may assume that $d(M)$ is so large that $F'$ is incompressible in $M$.  This means that we have a connected submanifold $M_2\cup M_F$ of $E\times I$ bounded by the incompressible surface $F'$.  As each component of a closed incompressible surface in $E\times I$ is parallel to $E$, the connected submanifold $M_2\cup M_F$ must be isotopic to $E\times I$ in $E\times I$.  Thus after isotopy, we may assume $M_2\cup M_F=E\times I$ and $F'=E\times\partial I$.

The compression body $M_F$ can be obtained by adding 1--handles to a small product neighborhood of $F'=E\times\partial I$.  So there is a graph $G$ properly embedded in $E\times I$ which corresponds to the 1--handles in $M_F$, such that after isotopy $M_F=\overline{N(G\cup(E\times\partial I))}$, where $N(G\cup(E\times\partial I))$ is a regular neighborhood of $G\cup(E\times\partial I)$ in $E\times I$.  Hence we may view $M_2=(E\times I)-N(G\cup(E\times\partial I))$.  Since the compression body $M_F$ is connected, the graph $G$ connects the two components of $F'=E\times\partial I$.

Our goal is to use the intersection of $M_2$ with a vertical annulus in $E\times I$ to construct an essential surface in $M_2$, and then apply Lemma~\ref{Ldist}.  We first show that a vertical annulus in $E\times I$ cannot be totally isotoped into $M_F=\overline{N(G\cup(E\times\partial I))}$.

We claim that $M_F=\overline{N(G\cup(E\times\partial I))}$ does not contain a properly embedded incompressible annulus $A$ whose two boundary circles lie in different components of $F'=E\times\partial I$.  As $M_F=\overline{N(G\cup(E\times\partial I))}$, after some handle slides if necessary, we may assume that there is a point $x$ in the graph $G$ that separates $E\times\partial I$ in the sense that no component of $G-x$ connects the two components of $F'=E\times\partial I$.  This means that there is a compressing disk $D$ for $F$ in $M_F$ such that $D$ is separating in $M_F$ and the two components of $F'$ lie in different components of $M_F-D$.  Suppose there is a properly embedded annulus $A$ described above.  As $\partial D\subset F$ and $\partial A\subset F'$, $A\cap\partial D=\emptyset$.  This means that $A\cap D$ (if not empty) consists of simple closed curves.  Since $A$ is incompressible, any curve in $A\cap D$ must bound disks in both $A$ and $D$. Hence after some isotopies removing closed curves in $A\cap D$ that are trivial in both $A$ and $D$, we have $A\cap D=\emptyset$.  However, this is impossible since  $A$ connects the two components of $F'$ but the two components of $F'$ lie in different components of $M_F-D$. 

Let $A$ be an essential vertical annulus in $E\times I$.  We may assume either $A\cap G=\emptyset$ or $A\cap G$ consists of a finite number of points in $\Int(A)$.  Hence $P=\overline{A-M_F}$ is a planar surface properly embedded in $M_2$.  After some standard cutting and pasting as in the proof of Corollary~\ref{Cirred}, we may assume $P$ is incompressible in $M_2$.  

The conclusion earlier says $M_F$ does not contain a properly embedded incompressible annulus whose two boundary circles lie in different components of $F'=E\times\partial I$.  So $A$ cannot be isotoped totally into $M_F$ and we cannot push $P$ into $M_F$.  This means that, after $\partial$-compressions on $P$, we obtain an essential planar surface $Q$ ($Q\ne\emptyset$) properly embedded in $M_2$. Since we can view a $\partial$-compression on $P$ as part of an isotopy on $A$ pushing the $\partial$-compressing disk into $M_F$, we may view $Q$ as a possibly disconnected subsurface of $A$ and $Q=A\cap M_2$. 

Next we show that there is a curve $\gamma_Q\subset\partial Q$ such that $d(\gamma_Q,\mathcal{D}_1)\le 1$. 
Since $F$ is incompressible in $M_2$, no component of $Q$ is  a disk.  If a component of $Q$ is not an essential subannulus of $A$, then there is a component $\gamma_Q$ of $\partial Q$ that  bounds a disk in $A$.  By Proposition~\ref{Ptrivial}, $d(\gamma_Q, \mathcal{D}_1)\le 1$.  If every component of $Q$ is an essential subannulus of $A$, then every component of $A-\Int(Q)$ is also an essential subannulus of $A$.  Let $A'$ be a component of $A-\Int(Q)$ that contains a boundary circle of $A$.  So $A'$ is an annulus properly embedded in $M_F$ with one component of $\partial A'$ in $F'=E\times\partial I$ and the other component of $\partial A'$, denoted by $\gamma_Q$, in $\partial Q\subset F$.  Since $\partial A$ is essential, $A'$ is incompressible in $M_F$.  If $A'$ is $\partial$-compressible in $M_F$, then  an essential arc of $A'$ bounds a $\partial$-compressing disk and hence has both endpoints in the same boundary component of $M_F$, which implies that $\partial A'$ lies in the same boundary component of $M_F$.  However, by our choice of $A'$, one boundary circle of $A'$ lies in $F'=E\times\partial I$ and the other boundary circle of $A'$ lies in $F$.  So $A'$ must be $\partial$-incompressible.  Hence $A'$ is an essential annulus in $M_F$.  After some standard cutting and pasting, one can always find a compressing disk for $F$ in $M_F$ disjoint from any essential annulus in $M_F$.  Thus $d(\gamma_Q,\mathcal{D}_1)\le 1$.  Hence, in any case, there is a curve $\gamma_Q\subset\partial Q$ such that $d(\gamma_Q,\mathcal{D}_1)\le 1$.

Note that $M_2$ cannot be a twisted $I$--bundle, since $E\times I$ does not contain any closed embedded non-orientable surface.  Let $\Omega_2$ be the essential surface used in defining $d(M)$.  Since $Q$ is a planar surface, by Lemma~\ref{Ldist}, there is a number $K'$ depending on $g(\Omega_2)$ such that $d(\partial\Omega_2\cap F,\gamma_Q)\le K'$.  Therefore 
$d(M)=d(\partial\Omega_2\cap F,\mathcal{D}_1)\le d(\partial\Omega_2\cap F,\gamma_Q)+d(\gamma_Q,\mathcal{D}_1)\le K'+1$.
\end{proof}

\begin{lemma}\label{Lproduct}
Let $E$ be a closed orientable incompressible surface of genus $g$ in $M$.  Then there is a number $K$ depending only on $g$ and $M_2$ such that if $F$ lies in a product neighborhood of $E$ in $M$, then $d(M)<K$.
\end{lemma}
\begin{proof} We may suppose $E\subset\Int(M)$ ($E$ may be parallel to a boundary component of $M$). 
Let $E\times I$ be a product neighborhood of $E$ in $\Int(M)$ and suppose $F\subset\Int(E\times I)$.  By Lemma~\ref{Lknotted}, we may assume $M_2\not\subset E\times I$.  So at least one component of $E\times\partial I$ lies in $M_2$. 

Since both $M_1$ and $M_2$ are irreducible and $M\not\cong S^3$, $F$ does not lie in a 3--ball in $E\times I$ and hence we can find a vertical annulus $A$ of $E\times I$ that cannot be isotoped disjoint from $F$. Let $N=M_2\cap (E\times I)$ and $P_2=A\cap N$.  Since at least one component of $E\times\partial I$ lies in $M_2$,  one or two components of $\partial P_2$ lie in $E\times\partial I$.

Let $D$ be a compressing disk for $F$ in $M_1$. 
Since $E$ is incompressible in $M$ and $F\subset\Int(E\times I)$, $D$ can be isotoped into $E\times I$.  By shrinking $D$ to be sufficiently small, we may assume $D\cap A=\emptyset$ and hence $\partial P_2\cap\partial D=\emptyset$.  This means that $d(\gamma,\mathcal{D}_1)\le 1$ for any component $\gamma$ of $\partial P_2\cap F=A\cap F$. Moreover, since $M$ is irreducible, after some standard cutting and pasting as in the proof of Corollary~\ref{Cirred}, we may assume that $P_2$ is incompressible in $N$.

Note that $\partial N$ consists of $F$ and one or both components of $E\times\partial I$.  Moreover, a component of $\partial P_2$ lies in $E\times\partial I$ and this implies that $P_2$ is not totally in a product neighborhood of $F$ in $N$.  
So we can apply Lemma~\ref{Lpartial}, setting $M'$, $P$ and $P_2$ in Lemma~\ref{Lpartial} to be $E\times I$, $A$ and $P_2$ above respectively.
 After performing some $\partial$-compressions on $P_2$ in $N$, we obtain an essential surface $Q$ such that $Q\cap F\ne\emptyset$ and $d(Q\cap F, (A\cap F)\cup\mathcal{D}_1)\le 3-\chi(A)=3$.  So there is a component $\delta$ of $Q\cap F$ such that $d(\delta,(A\cap F)\cup\mathcal{D}_1)\le 3$.  Since $d(\gamma,\mathcal{D}_1)\le 1$ for any component $\gamma$ of $\partial P_2\cap F=A\cap F$,  $d(\delta,\mathcal{D}_1)\le 4$ for some component $\delta$ of $Q\cap F$.

If $M_2$ is a twisted $I$--bundle over a closed non-orientable surface, since a component of $E\times\partial I$ lies in $M_2$ and is incompressible in $M_2$, $E$ must be parallel to $\partial M_2=F$.  However, this contradicts that $F$ is compressible in $M_1$ but $E$ is incompressible in $M$.  Thus  $M_2$ cannot be a twisted $I$--bundle.

Let $\Omega_2$ be the surface in $M_2$ used in defining $d(M)=d(\mathcal{D}_1, \partial\Omega_2\cap F)$.  Note that $\partial N$ is incompressible in $M_2$, since $\partial N$ consists of $F$ and one or both components of $E\times\partial I$.  So we may assume $\Omega_2\cap N$ is an essential subsurface of $\Omega_2$ and $-\chi(\Omega_2\cap N)\le -\chi(\Omega_2)$.   $N$ cannot be an $I$--bundle, since $F\subset\partial N$ is compressible in $M$ but $\partial N-F\subset E\times\partial I$ is incompressible in $M$.  By our construction of $P_2$ and $Q$, $Q$ is a planar surface in $N$ with all but one or two boundary components in $F$.  Thus by Lemma~\ref{Ldist}, $d(\partial\Omega_2\cap F, Q\cap F)\le K'$ for some $K'$ depending only on $\chi(\Omega_2)$.  Since $d(\delta,\mathcal{D}_1)\le 4$ for some component $\delta$ of $Q\cap F$, $d(M)=d(\partial\Omega_2\cap F,\mathcal{D}_1)\le d(\partial\Omega_2\cap F,\delta)+d(\delta,\mathcal{D}_1)\le K'+5$.
\end{proof}

Let $M=N_0\cup_{F_1}N_1\cup_{F_2}\dots\cup_{F_m}N_m$ be the decomposition in Theorem~\ref{TST} given by the untelescoping of an irreducible Heegaard splitting.  Recall that by Lemma~\ref{Llarge}, we have assumed that $F\subset\Int(N_j)$ for some $j$.  In the next two lemmas, we discuss the case that $F$ is also disjoint from the Heegaard surface $P_j$ of $N_j$.

\begin{lemma}\label{LM2} 
Let $N_j$ be the submanifold of $M$ between $F_j$ and $F_{j+1}$ in the untelescoping construction as above (we assume $N_j$ is connected), and let $P_j$ be the strongly irreducible Heegaard surface of $N_j$.  Suppose $F\subset\Int(N_j)$ and $F\cap P_j=\emptyset$.  Then there is a number $K$ such that if $d(M)>K$, $P_j\subset\Int(M_2)$.
\end{lemma}
\begin{proof}
Since $P_j\cap F=\emptyset$, $P_j$ lies in either $\Int(M_1)$ or $\Int(M_2)$. Suppose the lemma is not true and $P_j\subset\Int(M_1)$.  Let $X$ and $Y$ be the two compression bodies in the Heegaard splitting of $N_j$ along $P_j$.  We may suppose $F\subset\Int(X)$ and let $Z=X\cap M_1$.  Since $F\subset\Int(X)$ and $P_j\subset\Int(M_1)$, $F$ and $P_j$ are boundary components of $Z$. So we may view $Z\cup M_2$ as a submanifold of $M$.  By our construction, $X\subset Z\cup M_2$.   Since $P_j$ is compressible on both sides, there is a compressing disk $D$ for $P_j$ in $X\subset Z\cup M_2$.  We claim that $D$ can be isotoped into $Z$ if $d(M)$ is large.

Suppose on the contrary that $D$ cannot be isotoped into $Z$.  So $D\cap F\ne\emptyset$ after any isotopy on $D$ and we may assume $|D\cap F|$ is minimal in the isotopy class of $D$.  Let $Q$ be a component of $D\cap M_2$.  Since $D\cap F\ne\emptyset$ and $|D\cap F|$ is minimal, we may assume $Q$ cannot be pushed into $M_1$ and  $Q$ is incompressible in $M_2$ after isotopy.  As in the proofs of Corollary~\ref{Cirred} and Corollary~\ref{Cincomp}, we can perform some $\partial$--compressions on $Q$ in $M_2$ and obtain an essential planar surface $Q'$ in $M_2$.  We may regard $Q'$ as a subsurface of $D$.  Since every component of $\partial Q'$ bounds a disk in $D$, by Proposition~\ref{Ptrivial}, for any component $\gamma$ of $\partial Q'$, $d(\gamma,\mathcal{D}_1)\le 1$.  Now similar to the proofs of Corollary~\ref{Cirred} and Corollary~\ref{Cincomp}, this implies that $d(M)\le K$ for some $K$ depending only on $M_2$.  Thus if $d(M)$ is sufficiently large, every compressing disk of $P_j$ in $Z\cup M_2$ can be isotoped into $Z$.

Let $W$ be the surface obtained by maximally compressing $P_j$ in $Z$ and removing all resulting 2--sphere components.  
Since a maximal compression on $P_j$ in $X$ yields $\partial_-X$, the conclusion above implies that 
 $W$ is parallel to $\partial_-X$.  This means that $F$ must lie in a product region bounded by $W$ and $\partial_-X$.  Now Lemma~\ref{LM2} follows from Lemma~\ref{Lproduct}.
\end{proof}

\begin{lemma}\label{Lboth} 
Let $N_j$ be the submanifold of $M$ between $F_j$ and $F_{j+1}$ in the untelescoping construction as above (we assume $N_j$ is connected), and let $P_j$ be the strongly irreducible Heegaard surface of $N_j$. 
Suppose $F\subset\Int(N_j)$, $P_j\subset\Int(M_2)$ and $P_j$ is compressible on both sides in $M_2$.  Then there is a number $K$ such that, if $d(M)>K$, $F$ is isotopic to a middle surface of a compression body in the Heegaard splitting along $P_j$ (see Definition~\ref{Dmiddle}).
\end{lemma}
\begin{proof} 
Let $X$ and $Y$ be the two compression bodies in the Heegaard splitting of $N_j$ along $P_j$.  As $P_j\subset\Int(M_2)$, $F\cap P_j=\emptyset$.  Since $F\subset\Int(N_j)$ and $F\cap P_j=\emptyset$, we may suppose $F\subset\Int(X)$.
Let $Z=X\cap M_2$.  So $P_j$ and $F$ are boundary components of $Z$. Since $P_j$ is compressible on both sides in $M_2$, $P_j$ is compressible in $Z$. 
Let $P'$ be the surface obtained by maximally compressing $P_j$ in $Z$ and removing all resulting 2--sphere components.  So $P'$ is incompressible in $Z$.  Since $P_j$ is strongly irreducible, Casson-Gordon \cite{CG} implies that $P'$ is also incompressible on the other side.  Hence $P'$ is incompressible in $M_2$.

Let $N$ be the submanifold of $Z$ between $F$ and $P'$.  So $\partial N$ is incompressible in $M_2$ If $N$ is an $I$--bundle (i.e., if $F$ is parallel to a component of $P'$), then by our construction, $F$ is a middle surface of the compression body $X$ and the lemma holds.  Next we suppose $N$ is not an $I$--bundle. 

If $P'$ is parallel to $\partial_-X$, then by our construction, $F$ lies in the product region bounded by $P'$ and $\partial_-X$.  By Lemma~\ref{Lproduct}, we may suppose $d(M)$ is so large that $F$ does not lie in a product neighborhood of $\partial_-X$.  Thus we may assume $P'$ is not parallel to $\partial_-X$, and this means that $P'$ must be compressible in $X$ but incompressible in $Z$.  Let $D$ be a compressing disk for $P'$ in $X$.  By the construction of $P'$, $D\cap F\ne\emptyset$ after any isotopy on $D$.  Let $Q$ be the component of $D\cap N$ that contains $\partial D$.  After isotopy, we may assume $Q$ is an essential surface in $N$.  Note that one component of $\partial Q$ (i.e., $\partial D$) lies in $P'$ and all other components of $\partial Q$ lie in $F$.  Moreover, $Q\cap F=\partial Q-\partial D\ne\emptyset$ and every curve of $Q\cap F$ bounds a subdisk of $D$. By Proposition~\ref{Ptrivial}, $d(\mathcal{D}_1,\gamma)\le 1$ for every component $\gamma$ of $Q\cap F$.

Since $P'$ is incompressible in $M_2$ and $F$ is not parallel to $P'$, $M_2$ cannot be a twisted $I$--bundle. 
Let $\Omega_2$ be the fixed essential surface in $M_2$ used in defining $d(M)=d(\mathcal{D}_1,\partial\Omega_2\cap F)$.  Since $P'$ is incompressible in $M_2$, we may assume $\Omega'=\Omega_2\cap N$ is an essential subsurface of $\Omega_2$ and $\Omega'$ is essential in $N$.  Thus $-\chi(\Omega')\le-\chi(\Omega_2)$.

By Lemma~\ref{Ldist}, there is a number $K'$ depending on $\chi(\Omega_2)$  such that $d(Q\cap F, \partial\Omega'\cap F)\le K'$.  Let $\gamma$ be a component of $Q\cap F$ realizing $d(\gamma, \partial\Omega'\cap F)=d(Q\cap F, \partial\Omega'\cap F)$.  Since $\partial\Omega'\cap F=\partial\Omega_2\cap F$ and $d(\mathcal{D}_1,\gamma)\le 1$ for every $\gamma$ in $Q\cap F$,  $d(M)=d(\mathcal{D}_1,\partial\Omega_2\cap F)\le d(\mathcal{D}_1,\gamma)+d(\gamma,\partial\Omega_2\cap F)\le 1+K'$.
\end{proof}

Lemma~\ref{LM2} and Lemma~\ref{Lboth} say that if $F\subset N_j$ and $F\cap P_j=\emptyset$, then either Theorem~\ref{Tmain} holds, or (1) $P_j$ lies in $M_2$ and (2) $P_j$ cannot be compressible on both sides in $M_2$.

\section{Intersection of $F$ with sweepout surfaces}\label{Slast}
Let $M=M_1\cup_F M_2$ be as in section~\ref{Scase3} and $S$ an unstabilized genus $g$ Heegaard surface of $M$.  As in Theorem~\ref{TST}, let $M=N_0\cup_{F_1}N_1\cup_{F_2}\dots\cup_{F_m}N_m$ be the decomposition given by the untelescoping of $S$, where each $F_i$ is incompressible in $M$. 
By Lemma~\ref{Llarge}, we may assume that $d(M)$ is so large that each $F_i$ is disjoint from $F$.  Suppose $F\subset \Int(N_j)$.  Without loss of generality, we may assume $N_j$ is connected.  Let $P_i=P$ be the strongly irreducible Heegaard surface of $N_i$ in the untelescoping.

Let $F'$ be the surface obtained by maximally compressing $F$ in $M_1$ and deleting all the resulting 2--sphere components. 
We consider the compression body $M_F$ bounded by $F'$ and $F$.  So $\partial_+M_F=F$ and $\partial_-M_F=F'$.  Since $F\subset\Int(N_j)$ and $\partial N_j$ is incompressible in $M$, every compressing disk of $F$ in $M_1$ can be isotoped into $N_j$.  Thus we may assume $M_F\subset\Int(N_j)$.

Let $X$ and $Y$ be the two compression bodies in the Heegaard splitting of $N_j$ along $P_j=P$.  Let graphs $G_X$ and $G_Y$ be the cores of the compression bodies $X$ and $Y$ respectively, $\Sigma_X=G_X\cup\partial_-X$ and $\Sigma_Y=G_Y\cup\partial_-Y$, such that $N_j-(\Sigma_X\cup\Sigma_Y)\cong P\times(0,1)$.  We consider the sweepout $f:P\times I\to N_j$ such that $f|_{P\times(0,1)}$ is an embedding, $f(P\times\{0\})=\Sigma_X$ and $f(P\times\{1\})=\Sigma_Y$.  We denote $f(P\times\{t\})$ by $P^t$. Each $P^t$ is isotopic to $P_j=P$ if $t\in(0,1)$.

Similarly, let graph $G_F$ be the core of the compression body $M_F$ and $\Sigma_F=G_F\cup F'$ such that $M_F-\Sigma_F\cong F\times(0,1]$.  After isotopy, we may assume that the graphs $G_X$, $G_Y$ and $G_F$ are pairwise disjoint in $N_j$. 
We may suppose $F'$ is transverse to the graphs $G_X$ and $G_Y$ and suppose $F'$ is transverse to each $P^t$ except for finitely many levels $t$ where $F'\cap P^t$ contains exactly one center or saddle tangency.  We may suppose $G_F\cap P^t$ consists of finitely many points for each $t\in I$.  Moreover, we may suppose $M_F$ is a small neighborhood of $\Sigma_F=G_F\cup F'$ in $N_j$ and suppose $F$ is transverse to each $P^t$ except for finitely many levels $t$ where $F\cap P^t$ contains exactly one center or saddle tangency.  As $M_F$ is a small neighborhood of $\Sigma_F$, $F$ is transverse to the graphs $G_X$ and $G_Y$.

We use $\Lambda$ to denote the union of $\partial I=\{0,1\}$ and the levels $t\in I$ where $F$ is not transverse to $P^t$.

For any $P^t$, $t\in(0,1)$, we use $X_t$ (resp. $Y_t$) to denote
the closure of the component of $N_j-P^t$ that contains $\Sigma_X$ 
(resp. $\Sigma_Y$).  Recall that $P=P_j$ is a Heegaard surface of $N_j$ and $P^t$ is parallel to $P$.  So if $t\in(0,1)$, $X_t$ and $Y_t$ are the two compression bodies corresponding to $X$ and $Y$ respectively.

\begin{labelling}  A number $t\in I-\Lambda$ 
 is labelled $X$ (resp.~$Y$) if, there is an essential curve $\gamma$ in $P^t$ such that
\begin{enumerate}
\item $\gamma$ bounds a compressing disk in $X_t$ (resp.~$Y_t$) and
\item $\gamma\subset M_2$ and $\gamma$ bounds an embedded disk $D$ in $M_2$ that is transverse to $P^t$.
\end{enumerate}
This is a labelling for all $t\in I-\Lambda$ and we do not assign any label for $t$ if $F$ is tangent to $P^t$.  
\end{labelling}
\begin{remark}\label{Remark}
 The graph $G_F$ can be viewed as the core of the 1--handles attached to (a product neighborhood of) $F'=\partial_-M_F$ in the compression body $M_F$.  Recall that we have assumed $M_F$ is a small neighborhood of $F'\cup G_F$ and $F$ is a boundary component of this small neighborhood. 
For any point $t\notin\Lambda$, since we have assumed that $P^t$ intersects the graph $G_F$ in finitely many points and since $M_F$ is a small neighborhood of $F'\cup G_F$, 
one can always find a compressing disk (which can be chosen to be a cocore of the 1--handles corresponding to some point in $G_F-P^t$) for $F$ in $M_F$ that is disjoint from $P^t$. 
Thus $d(\gamma,\mathcal{D}_1)\le 1$ for any component $\gamma$ of $F\cap P^t$ that is essential in $F$.
\end{remark}

Theorem~\ref{Tmain} follows from the following 4 claims.

\begin{claim}\label{Claim1} There is a number $K_1$ such that if $d(M)>K_1$, then 
for a sufficiently small $\epsilon>0$, $\epsilon$ is labelled $X$ and $1-\epsilon$ is labelled $Y$.
\end{claim}
\begin{proof}
Suppose $d(M)$ is larger than the constant $K$ in Lemma~\ref{LM2}.
 
Since the compression body $M_F$ lies in $N_j$, $F$ is disjoint from $\partial N_j=\partial_-X\cup\partial_-Y$. 
The graph $G_X$ cannot totally lie in $M_1$, because otherwise $P^\epsilon$ lies in $M_1$ for a sufficiently small $\epsilon$, contradicting Lemma~\ref{LM2}.  Thus $G_X\cap\Int(M_2)\ne\emptyset$ and $X_\epsilon$ must have a compressing disk $D$ lying in $M_2$. Hence $\epsilon$ is labelled $X$.

We can apply the same argument to $G_Y$ and conclude that $1-\epsilon$ is labelled $Y$ for sufficiently small $\epsilon$.
\end{proof}

\begin{claim}\label{Claimno}
There is a number $K_2$ such that if $d(M)>K_2$, then every $t\in  I-\Lambda$ has a label $X$ or $Y$.
\end{claim}
\begin{proof}
Suppose on the contrary that some $t\in  I-\Lambda$ has no label.  We assume $d(M)$ is larger than the constant $K$ in Lemma~\ref{LM2}.

Let $P_2=P^t\cap M_2$. By Lemma~\ref{LM2}, $P_2\ne\emptyset$.  Our goal is to use $P_2$ to construct an incompressible surface in $M_2$ and then apply the inequalities in Lemma~\ref{Lpartial} and Lemma~\ref{Ldist} to get a bound on the distance $d(M)$.

We first suppose $P_2$ is compressible in $M_2$ and let $D$ be a compressing disk for $P_2$ in $M_2$.  Since $t$ is not labelled, $D$ cannot be a compressing disk for $P^t$ and $\partial D$ must be trivial in $P^t$ but essential in $P_2$.  We compress $P^t$ along $D$ and delete the resulting 2--sphere component.  Let $P'$ be the remaining surface after this operation.  Since $M$ is irreducible, $P'$ is isotopic to $P^t$.  Suppose $P'\cap M_2$ is still compressible in $M_2$ and let $D'$ be a compressing disk of $P'\cap M_2$ in $M_2$.  Suppose $\partial D'$ is essential in $P'$.  Since $D'\cap P'=\partial D'$, by the operation above and after a slight perturbation on $D'$ if necessary, we may view $\partial D'$ as an essential curve in $P^t$ bounding an embedded disk $D'\subset M_2$.  Since $D'$ may intersect the 2--sphere component that we eliminated in the operation above, $\Int(D')\cap P^t$ may not be empty.  Nonetheless, since $\partial D'$ is essential in $P^t$, by Scharlemann's no-nesting lemma \cite[Lemma 2.2]{S}, $\partial D'$ bounds a compressing disk for $P^t$ in $X_t$ or $Y_t$.  This means that $t$ is labelled $X$ or $Y$, contradicting our hypothesis.  Thus $\partial D'$ must also be trivial in $P'$ and we can perform the same operation on $P'$, i.e.~compress $P'$ along $D'$ and remove the resulting 2--sphere component.

After finitely many such operations, we may assume that $P_2=P^t\cap M_2$ is incompressible in $M_2$.  If $P^t\cap F=\emptyset$, then by Lemma~\ref{LM2}, $P^t\subset\Int(M_2)$ and $P_2=P^t$.  However, since $P^t$ is separating in $N_j$ and compressible on both sides, $P_2=P^t$ must be compressible in $M_2$, a contradiction to our assumption on $P_2$.  Thus $P^t\cap F\ne\emptyset$.

By Lemma~\ref{LM2}, $P_2$ does not lie in a product neighborhood of $F$ (otherwise $P_2$ and hence $P^t$ can be isotoped into $M_1$). 
So after some $\partial$--compressions on $P_2$, we get an essential surface $Q$ properly embedded in $M_2$.  Now we apply Lemma~\ref{Lpartial}, setting $M'$, $P$ and $P_2$ in Lemma~\ref{Lpartial} to be $N_j$, $P^t$ and $P_2$ above respectively, and get $d(\partial Q,(P^t\cap F)\cup\mathcal{D}_1)\le 3-\chi(P^t)\le 2g+1$.  By Remark~\ref{Remark}, $d(\alpha,\mathcal{D}_1)\le 1$ for every component $\alpha$ of $P^t\cap F$.  Hence, $d(\partial Q, \mathcal{D}_1)\le  2g+2$. 

If $M_2$ is a twisted $I$--bundle, then $Q$ must be a vertical annulus.  Hence $d(M)\le d(\partial Q,\mathcal{D}_1)\le 2g+2$.

If $M_2$ is not a twisted $I$--bundle, let $\Omega_2$ be the fixed essential surface used in defining $d(M)$. 
As the genus of $Q$ is at most $g$, by Lemma~\ref{Ldist}, there is a $K'$ depending on $\Omega_2$ and $g$, such that $d(\partial\Omega_2\cap F,\partial Q)\le K'$.  Since $d(\partial Q,\mathcal{D}_1)\le 2g+2$, this means that $d(M)=d(\Omega_2\cap F,\mathcal{D}_1)\le d(\partial\Omega_2\cap F,\partial Q)+1+d(\partial Q,\mathcal{D}_1)\le K'+2g+3$.

Therefore if $d(M)$ is sufficiently large, every $t\in  I-\Lambda$ has a label.
\end{proof}

\begin{claim}\label{Claimboth}
If some $t\in  I-\Lambda$ is labelled both $X$ and $Y$, then Theorem~\ref{Tmain} holds.
\end{claim}
\begin{proof}
In this claim, we assume $d(M)$ is larger than the constants $K$ in Lemma~\ref{LM2} and Lemma~\ref{Lboth}.

Let $D$ be an embedded disk in $M_2$ transverse to $P^t$ with $\partial D\subset P^t\cap M_2$.  We call $D$ an \emph{almost compressing disk} for $X_t$ (resp.~$Y_t$) if $\partial D$ bounds a compressing disk in $X_t$ (resp.~$Y_t$).  

Suppose $t\in  I-\Lambda$ is labelled both $X$ and $Y$.  Then by definition, $M_2$ contains almost compressing disks $D_X$ and $D_Y$ for $X_t$ and $Y_t$ respectively.  Since $P^t$ is strongly irreducible, $\partial D_X\cap\partial D_Y\ne\emptyset$. 

Let $P_2=P^t\cap M_2$.  Similar to the proof of Claim~\ref{Claimno}, our goal is to use $P_2$ to construct either an incompressible or a strongly irreducible surface in $M_2$, and then apply the inequalities in Lemma~\ref{Lpartial} and Lemma~\ref{Ldist} to get a bound on the distance $d(M)$.  Although $P^t$ is strongly irreducible in $M$, $P_2$ may not be strongly irreducible in $M_2$ because the boundary curve of a compressing disk for $P_2$ in $M_2$ may be a trivial curve in $P^t$.  So we need to perform some operations on $P_2$ first.

Let $\Delta$ be a compressing disk for $P_2$ in $M_2$.  We say $\Delta$ is a \emph{trivial compressing disk} if $\partial\Delta$ is essential in $P_2$  but trivial in $P^t$.  Suppose a trivial compressing disk $\Delta$ lies in $X_t\cap M_2$ and there is an almost compressing disk $D_Y$ for $Y_t$ such that $\partial D_Y\cap\partial\Delta=\emptyset$.  Then we can compress $P^t$ along $\Delta$ and delete the resulting 2--sphere component.  As in Claim~\ref{Claimno}, the remaining surface $P'$ is isotopic to $P^t$.  Since $\partial D_Y\cap\partial\Delta=\emptyset$, $\partial D_Y\subset P'$ and $D_Y$ remains an almost compressing disk for $P'$.  

For any almost compressing disk $D_X$ for $X_t$, if a component $\gamma$ of $\Int(D_X)\cap P^t$ is essential in $P^t$, then by Scharlemann's no-nesting lemma \cite[Lemma 2.2]{S}, $\gamma$ must bound a compressing disk for $P^t$. Since $P^t$ is strongly irreducible and $\partial D_X$ bounds a compressing disk in $X_t$, the subdisk of $D_X$ bounded by $\gamma$ must also be an almost compressing disk for $X_t$.  Thus we may choose an almost compressing disk $D_X$ for $X_t$ so that every component of $\Int(D_X)\cap P^t$ is trivial in $P^t$.  Since $\partial D_X$ bounds a compressing disk in $X_t$, this implies that a small neighborhood of $\partial D_X$ in $D_X$ lies in $X_t$.  Now we consider $D_X\cap\Delta$, where $\Delta$ is the trivial compressing disk in $X_t\cap M_2$ above.  If $D_X\cap\Delta\ne\emptyset$, similar to the proof of Lemma~\ref{Lside}, we can push the arcs in $D_X\cap\Delta$ across $\Delta$.  More specifically, we may suppose $D_X\cap\Delta$ does not contain any closed curve and let $\alpha$ be an arc in $D_X\cap\Delta$ that is outermost in $\Delta$.  Then $\alpha$ and a subarc of $\partial\Delta$ bound a subdisk $E$ of $\Delta$ and $\Int(E)\cap D_X=\emptyset$.  Since a small neighborhood of $\partial D_X$ in $D_X$ lies in $X_t$ and $\Delta\subset X_t$, we have $E\cap D_X=\alpha$.  So, similar to the proof of Lemma~\ref{Lside}, we can perform an isotopy by pushing $\alpha$ and $D_X$ across $E$ to eliminate $\alpha$.  After the operation pushing $D_X$ across $E$ above, $D_X$ becomes either one or two disks depending on whether or not both endpoints of $\alpha$ lie in $\partial D_X$.  Since $\partial D_X$ is essential in $P^t$ and each component of $\Int(D_X)\cap P^t$ is trivial in $P^t$, after the operation, the boundary curve of at least one resulting disk is essential in $P^t$.  Hence after the operation pushing $D_X$ across $E$ above, we obtain a new almost compressing disk for $X_t$  with fewer intersection arcs with $\Delta$. 
After finitely many these operations, we can construct an almost compressing disk $D_X'$ (for $X_t$) that is disjoint from $\Delta$.  

The arguments above say that if there is a trivial compressing disk $\Delta$ in $X_t\cap M_2$ such that $\partial\Delta\cap\partial D_Y=\emptyset$ for some almost compressing disk $D_Y$ for $Y_t$, then after compressing $P^t$ along $\Delta$ and deleting the 2--sphere component, the resulting surface still has two almost compressing disks for $X_t$ and $Y_t$ respectively.  
Therefore, after finitely many such operations on trivial compressing disks as above, we may assume that for each trivial compressing disk $\Delta$, if $\Delta\subset X_t$ then $\partial\Delta\cap\partial D_Y\ne\emptyset$ for every almost compressing disk $D_Y$ for $Y_t$, and if $\Delta\subset Y_t$ then $\partial\Delta\cap\partial D_X\ne\emptyset$ for every almost compressing disk $D_X$ for $X_t$.  Note that this implies that every curve of $P^t\cap F$ must be essential in $F$, because otherwise the subdisk of $F$ bounded by an innermost such curve is either a trivial compressing disk disjoint from all almost compressing disks, or a compressing disk of $X_t$ (resp.~$Y_t$) disjoint from an almost compressing disk $D_Y$ (resp.~$D_X$), which contradicts that $P^t$ is strongly irreducible.

Next we show that $P_2=P^t\cap M_2$ has compressing disks in both $X_t\cap M_2$ and $Y_t\cap M_2$.  Suppose $P_2$ does not have any compressing disk lying in $X_t\cap M_2$.  Let $D_X$ be an almost compressing disk for $X_t$ and we may assume $|\Int(D_X)\cap P^t|$ is minimal among all almost compressing disks for $X_t$.  If $\Int(D_X)\cap P^t=\emptyset$, then $D_X$ is a compressing disk for $P_2$ lying in $X_t\cap M_2$, contradicting our assumption.  So we may suppose $\Int(D_X)\cap P^t\ne\emptyset$.  Let $\gamma$ be an innermost component of $\Int(D_X)\cap P^t$ and let $d_\gamma$ be the subdisk of $D_X$ bounded by $\gamma$.  If $\gamma$ is trivial in $P_2$, then we can perform a simple isotopy on $D_X$ to remove $\gamma$ and get a contradiction to the minimality assumption of $|\Int(D_X)\cap P^t|$.  Thus $\gamma$ is essential in $P_2$ and $d_\gamma$ is a compressing disk for $P_2$.  Since 
we have assumed that $P_2$ does not have any compressing disk lying in $X_t\cap M_2$,  $d_\gamma\subset Y_t\cap M_2$.  If $\gamma$ is also essential in $P^t$, then $d_\gamma$ is a compressing disk for $P^t$ in $Y_t$.  However, since $\partial D_X$ bounds a compressing disk for $P^t$ in $X_t$ and $\gamma\cap\partial D_X=\emptyset$, this contradicts that $P^t$ is strongly irreducible.  Hence $\gamma$ must be trivial in $P^t$ and $d_\gamma$ is a trivial compressing disk for $P_2$ in $Y_t$, but this contradicts our earlier assumption that every trivial compressing disk in $Y_t$ intersects every almost compressing disk for $X_t$ because $\gamma\cap\partial D_X=\emptyset$.  Therefore, $P_2=P^t\cap M_2$ must have compressing disks in both $X_t\cap M_2$ and $Y_t\cap M_2$.

Suppose $P_2$ is not strongly irreducible in $M_2$.  Then there are compressing disks $\Delta_X$ and $\Delta_Y$ for $P_2$ in $M_2$ such that $\Delta_X\subset X_t$ and $\Delta_Y\subset Y_t$ and $\partial\Delta_X\cap\partial\Delta_Y=\emptyset$.  By our assumptions above, both $\Delta_X$ and $\Delta_Y$ must be trivial compressing disks.  Now we compress $P^t$ along $\Delta_X$ and $\Delta_Y$ simultaneously and delete the two resulting 2--sphere components.  The remaining surface $P'$ is isotopic to $P^t$.  Suppose $P'\cap M_2$ has an almost compressing disk $D'$.  As in Claim~\ref{Claimno}, after some perturbation, we may view $\partial D'$ as an essential curve in $P^t$ and view $D'$ as an almost compressing disk of $P_2$.  However, since $D'\cap\Delta_X=\emptyset$ and $D'\cap\Delta_Y=\emptyset$ after isotopy, this contradicts our earlier assumption that every trivial compressing disk must intersect every almost compressing disk on the other side.  So $P'$ does not have any almost compressing disk in $M_2$, and this implies that every compressing disk of $P'\cap M_2$ in $M_2$ is a trivial compressing disk for $P'$.  We can compress $P'$ along each trivial compressing disk of $P'\cap M_2$ in $M_2$ and delete the resulting 2--sphere component.  By the argument above, the resulting surface does not have any almost compressing disk in $M_2$ neither.  Therefore, after finitely many such operations, we obtain a surface $P''$ isotopic to $P^t$ and $P''\cap M_2$ is incompressible in $M_2$.

The arguments above imply that, after some isotopies/operations on $P^t$ described above, we may assume that  $P_2=P^t\cap M_2$ is either strongly irreducible or incompressible in $M_2$.  If $P^t\cap F=\emptyset$ after the operations above, then by Lemma~\ref{LM2} $P^t\subset\Int(M_2)$ and hence $P_2=P^t$.  Since $P^t$ is separating in $N_j$ and compressible on both sides, $P^t\cap F=\emptyset$ implies that $P_2=P^t$ cannot be incompressible in $M_2$.  Hence $P_2=P^t$ is strongly irreducible and in particular $P_2$ is compressible on both sides in $M_2$.  In this case, by Lemma~\ref{Lboth}, $F$ is isotopic to a middle surface of the compression body $X^t$ or $Y^t$ and Theorem~\ref{Tmain} holds.  

Suppose Theorem~\ref{Tmain} fails, then the argument above implies that $P^t\cap F\ne\emptyset$. 

By Lemma~\ref{LM2}, $P_2$ does not lie in a product neighborhood of $F$ (otherwise $P_2$ and hence $P^t$ can be isotoped into $M_1$).  As $P_2$ is either strongly irreducible or incompressible in $M_2$, Claim~\ref{Claimboth} basically follows from Lemma~\ref{Lpartial}. 
By Lemma~\ref{Lpartial} (setting $M'$, $P$ and $P_2$ in Lemma~\ref{Lpartial} to be $N_j$, $P^t$ and $P_2$ above respectively), we can perform some $\partial$-compressions on $P_2$ in $M_2$ and obtain a surface $Q$ which is either essential or strongly irreducible and $\partial$-strongly irreducible, such that
 $d(\partial Q,(P^t\cap F)\cup\mathcal{D}_1)\le 3-\chi(P^t)\le 2g+1$.  By  Remark~\ref{Remark}, $d(\alpha,\mathcal{D}_1)\le 1$ for every component $\alpha$ of $P^t\cap F$.  Hence, $d(\partial Q, \mathcal{D}_1)\le 2g+2$. 

Suppose $M_2$ is not a twisted $I$--bundle and 
let $\Omega_2$ be the essential surface used in defining $d(M)$.  As the genus $g(Q)\le g$, by Lemma~\ref{Ldist}, there is a number $K'$ depending on $\Omega_2$ and $g$ such that $d(\partial\Omega_2\cap F,\partial Q)\le K'$.  Thus $d(M)=d(\partial\Omega_2\cap F,\mathcal{D}_1)\le d(\partial\Omega_2\cap F,\partial Q)+1+d(\partial Q, \mathcal{D}_1)\le K'+2g+3$.

If $M_2$ is a twisted $I$--bundle, then we can apply Lemma~\ref{Lbundle} instead of Lemma~\ref{Ldist} in the argument above and get the same inequality.

Therefore, if $d(M)$ is sufficiently large and some $t\in  I-\Lambda$ is labelled both $X$ and $Y$, $F$ must be isotopic to a middle surface in $X$ or $Y$ as in Lemma~\ref{Lboth} and Theorem~\ref{Tmain} holds.
\end{proof}

\begin{claim}\label{Claim3}
Suppose every $t\in  I-\Lambda$ is labelled, then Theorem~\ref{Tmain} holds.
\end{claim}
\begin{proof}
By Claim~\ref{Claimboth}, we may assume that no $t$ is labelled both $X$ and $Y$.
By Claim~\ref{Claim1}, as $t$ increases from $\epsilon$ to $1-\epsilon$, its label changes from $X$ to $Y$.   As $t\in  I-\Lambda$ is labelled, then  there is a number $t_0\in \Lambda$ such that $t_0-\epsilon$ is labelled $X$ and $t_0+\epsilon$ is labelled $Y$ for sufficiently small $\epsilon>0$.   Since $t_0 \in\Lambda$, $F\cap P^{t_0}$ contains a single tangency.  Since $t_0-\epsilon$ and $t_0+\epsilon$ have different labels, the tangency in $F\cap P^{t_0}$ must be a saddle tangency.

Let $F\times J$ be a small product neighborhood of $F$ in $M$, where $J$ is a closed interval, and let $F^+$ and $F^-$ be the two components of $F\times\partial J$.  $F^+$ and $F^-$ are parallel and close to $F$ but lie on different sides of $F$.  By considering how the intersection curves change near a saddle tangency, it is easy to see that, we may choose $F^+$ and $F^-$ so that, for a sufficiently small $\epsilon$, the intersection patterns of $F^{\pm}\cap P^{t_0}$ and $F\cap P^{t_0\pm\epsilon}$ are the same.  In fact, we may assume $F^{+}\cup P^{t_0}$ and $F^{-}\cup P^{t_0}$ are isotopic to $F\cup P^{t_0+\epsilon}$ and $F\cup P^{t_0-\epsilon}$ in $M$ respectively.

Let $M_1^\pm$ and $M_2^\pm$ be components of the closure of $M-F^{\pm}$ corresponding to $M_1$ and $M_2$ respectively.  There are two subcases depending on whether $F^{\pm}$ lies in $M_1$ or $M_2$.

The first subcase is that $F^-\subset\Int(M_1)$ and $F^+\subset\Int(M_2)$.
In this subcase $M_1^-\subset\Int(M_1^+)$ and $M_2^+\subset\Int(M_2^-)$.  Since $t_0-\epsilon$ is labelled $X$ and since intersection pattern of $F^{-}\cap P^{t_0}$ is the same as $F\cap P^{t_0-\epsilon}$, there is a curve $\gamma_X$ in $P^{t_0}\cap M_2^-$ such that (1) $\gamma_X$ bounds a compressing disk in $X_{t_0}$, and (2) $\gamma_X$ bounds an almost compressing disk $D_X$ in $M_2^-$.  Similarly, since $t_0+\epsilon$ is labelled $Y$, there is a curve $\gamma_Y$ in $P^{t_0}\cap M_2^+$ such that (1) $\gamma_Y$ bounds a compressing disk in $Y_{t_0}$, and (2) $\gamma_Y$ bounds an almost compressing disk $D_Y$ in $M_2^+$.  Since  $M_2^+\subset\Int(M_2^-)$, $D_Y\subset M_2^+\subset M_2^-$.  So both $D_X$ and $D_Y$ lie in $M_2^-$.  As $F^{-}\cup P^{t_0}$ and $F\cup P^{t_0-\epsilon}$ are isotopic in $M$, there are a pair of almost compressing disks $D_X'$ and $D_Y'$ for $P^{t_0-\epsilon}$ in $M_2$ corresponding to $D_X$ and $D_Y$.  This means that $t_0-\epsilon$ is labelled both $X$ and $Y$.  Now Theorem~\ref{Tmain} follows from Claim~\ref{Claimboth}.

The second subcase is that $F^+\subset\Int(M_1)$ and $F^-\subset\Int(M_2)$. This subcase is basically the same as the first one.  One can simply interchange all the plus and minus signs and interchange all the labels $X$ and $Y$ in the proof above for the first subcase to conclude that $t_0+\epsilon$ is labelled both $X$ and $Y$.  By Claim~\ref{Claimboth}, Theorem~\ref{Tmain} holds.
\end{proof}

\end{document}